\documentclass[11pt,a4paper,leqno,twoside]{article}
\usepackage[a4paper,textwidth=13.5cm,textheight=19.8cm,centering]{geometry}
\usepackage[T1]{fontenc}
\usepackage{lmodern}
\usepackage{textcomp}
\usepackage{amsmath,amssymb,amsthm}
\usepackage{booktabs}
\usepackage{tikz}
\usepackage{float}
\usetikzlibrary{calc}

\newcommand{\rwportlabels}[2]{%
  \node[font=\scriptsize,anchor=south] at (#1-1.10,0.92) {$I_0^{#2}$};%
  \node[font=\scriptsize,anchor=south] at (#1-1.10,0.07) {$I_1^{#2}$};%
  \node[font=\scriptsize,anchor=south] at (#1-1.10,-0.78) {$I_2^{#2}$};%
  \node[font=\scriptsize,anchor=south] at (#1+1.10,0.92) {$O_0^{#2}$};%
  \node[font=\scriptsize,anchor=south] at (#1+1.10,0.07) {$O_1^{#2}$};%
  \node[font=\scriptsize,anchor=south] at (#1+1.10,-0.78) {$O_2^{#2}$};%
}

\newcommand{\rwlongwire}[2]{%
  \draw[#1] (#2+0.85,-0.85) -- (#2+1.60,-0.85) -- (#2+2.05,-1.25) -- (#2+5.25,-1.25) -- (#2+5.70,-0.85) -- (#2+6.45,-0.85);%
}
\newcommand{\rwlongwires}[1]{%
  \rwlongwire{#1}{0}%
  \rwlongwire{#1}{3.65}%
  \draw[#1] (8.15,-0.85) -- (8.90,-0.85) -- (9.35,-1.25) -- (12.55,-1.25) -- (13.00,-0.85) -- (13.20,-0.85);%
  \draw[#1] (11.80,-0.85) -- (12.55,-0.85) -- (13.00,-1.25) -- (13.20,-1.25);%
  \draw[#1] (-1.40,-0.85) -- (-0.85,-0.85);%
  \draw[#1] (-1.40,-1.25) -- (1.60,-1.25) -- (2.05,-0.85) -- (2.80,-0.85);%
}

\newcommand{\rwopaque}[2]{%
  \node[draw=black!65,fill=black!3,rounded corners=1.5pt,minimum width=1.30cm,
    minimum height=1.70cm,inner sep=0pt] at (#1,0) {};%
  \node[font=\small] at (#1,0) {$H^{#2}$};%
  \draw[black!65,line width=0.55pt] (#1-0.65,0.85) -- (#1-0.85,0.85);%
  \draw[black!65,line width=0.55pt] (#1-0.65,-0.85) -- (#1-0.85,-0.85);%
  \draw[black!65,line width=0.55pt] (#1+0.65,0.85) -- (#1+0.85,0.85);%
  \draw[black!65,line width=0.55pt] (#1+0.65,-0.85) -- (#1+0.85,-0.85);%
  \rwportlabels{#1}{#2}%
}

\newcommand{\ringopaquepanel}{%
\begin{tikzpicture}[x=1cm,y=1cm,line cap=round,line join=round]
  \tikzset{wire/.style={draw=black!65,line width=0.65pt}}
  \rwlongwires{wire}
  \foreach \x in {0,3.65,7.30}{%
    \draw[wire] (\x+0.85,0.85) -- (\x+2.80,0.85);
    \draw[wire] (\x+0.65,0) -- (\x+3.00,0);
  }
  \draw[wire] (-1.40,0.85) -- (-0.85,0.85);
  \draw[wire] (-1.40,0) -- (-0.65,0);
  \draw[wire] (11.80,0.85) -- (13.20,0.85);
  \draw[wire] (11.60,0) -- (13.20,0);
    \node[font=\small,anchor=west] at (-1.40,1.75) {(c)};
    \rwopaque{0}{i}
    \rwopaque{3.65}{i+1}
    \rwopaque{7.30}{i+2}
    \rwopaque{10.95}{i+3}
\end{tikzpicture}%
}

\newcommand{\ringhexpanel}{%
\begin{tikzpicture}[x=1cm,y=1cm,line cap=round,line join=round]
  \tikzset{wire/.style={draw=black!65,line width=0.65pt},
    hexedge/.style={draw=black!58,line width=0.9pt},
    border/.style={draw=black!30,line width=0.5pt,dash pattern=on 2.4pt off 1.6pt,rounded corners=1.5pt},
    dotn/.style={circle,fill=black,inner sep=1.1pt},
    lab/.style={font=\scriptsize,anchor=south},
    hlab/.style={font=\scriptsize,inner sep=0pt}}
  \draw[border] (-0.650,-1.020) rectangle (0.650,1.020);
  \draw[border] (3.000,-1.020) rectangle (4.300,1.020);
  \draw[border] (6.650,-1.020) rectangle (7.950,1.020);
  \draw[border] (10.300,-1.020) rectangle (11.600,1.020);
  \draw[hexedge] (-0.450,0.850) -- (0.450,0.850);
  \draw[hexedge] (0.450,0.850) -- (0.450,-0.850);
  \draw[hexedge] (0.450,-0.850) -- (0.180,-0.250);
  \draw[hexedge] (0.180,-0.250) -- (-0.180,-0.250);
  \draw[hexedge] (-0.180,-0.250) -- (-0.450,-0.850);
  \draw[hexedge] (-0.450,-0.850) -- (-0.450,0.850);
  \draw[hexedge] (3.200,0.850) -- (4.100,0.850);
  \draw[hexedge] (4.100,0.850) -- (4.100,-0.850);
  \draw[hexedge] (4.100,-0.850) -- (3.830,-0.250);
  \draw[hexedge] (3.830,-0.250) -- (3.470,-0.250);
  \draw[hexedge] (3.470,-0.250) -- (3.200,-0.850);
  \draw[hexedge] (3.200,-0.850) -- (3.200,0.850);
  \draw[hexedge] (6.850,0.850) -- (7.750,0.850);
  \draw[hexedge] (7.750,0.850) -- (7.750,-0.850);
  \draw[hexedge] (7.750,-0.850) -- (7.480,-0.250);
  \draw[hexedge] (7.480,-0.250) -- (7.120,-0.250);
  \draw[hexedge] (7.120,-0.250) -- (6.850,-0.850);
  \draw[hexedge] (6.850,-0.850) -- (6.850,0.850);
  \draw[hexedge] (10.500,0.850) -- (11.400,0.850);
  \draw[hexedge] (11.400,0.850) -- (11.400,-0.850);
  \draw[hexedge] (11.400,-0.850) -- (11.130,-0.250);
  \draw[hexedge] (11.130,-0.250) -- (10.770,-0.250);
  \draw[hexedge] (10.770,-0.250) -- (10.500,-0.850);
  \draw[hexedge] (10.500,-0.850) -- (10.500,0.850);
  \draw[wire] (0.450,0.850) -- (3.200,0.850);
  \draw[wire] (0.180,-0.250) -- (0.430,0.000) -- (3.220,0.000) -- (3.470,-0.250);
  \draw[wire] (4.100,0.850) -- (6.850,0.850);
  \draw[wire] (3.830,-0.250) -- (4.080,0.000) -- (6.870,0.000) -- (7.120,-0.250);
  \draw[wire] (7.750,0.850) -- (10.500,0.850);
  \draw[wire] (7.480,-0.250) -- (7.730,0.000) -- (10.520,0.000) -- (10.770,-0.250);
  \draw[wire] (0.450,-0.850) -- (1.600,-0.850) -- (2.050,-1.250) -- (5.250,-1.250) -- (5.700,-0.850) -- (6.850,-0.850);
  \draw[wire] (4.100,-0.850) -- (5.250,-0.850) -- (5.700,-1.250) -- (8.900,-1.250) -- (9.350,-0.850) -- (10.500,-0.850);
  \draw[wire] (-1.400,0.850) -- (-0.450,0.850);
  \draw[wire] (-1.400,0.000) -- (-0.430,0.000) -- (-0.180,-0.250);
  \draw[wire] (-1.400,-0.850) -- (-0.450,-0.850);
  \draw[wire] (-1.400,-1.250) -- (1.600,-1.250) -- (2.050,-0.850) -- (3.200,-0.850);
  \draw[wire] (11.400,0.850) -- (13.200,0.850);
  \draw[wire] (11.130,-0.250) -- (11.380,0.000) -- (13.200,0.000);
  \draw[wire] (11.400,-0.850) -- (12.550,-0.850) -- (13.000,-1.250) -- (13.200,-1.250);
  \draw[wire] (7.750,-0.850) -- (8.900,-0.850) -- (9.350,-1.250) -- (12.550,-1.250) -- (13.000,-0.850) -- (13.200,-0.850);
  \node[hlab] at (0.000,0.420) {$H^{i}$};
  \node[hlab] at (3.650,0.420) {$H^{i+1}$};
  \node[hlab] at (7.300,0.420) {$H^{i+2}$};
  \node[hlab] at (10.950,0.420) {$H^{i+3}$};
  \node[dotn] at (0.450,0.850) {};
  \node[dotn] at (0.450,-0.850) {};
  \node[dotn] at (0.180,-0.250) {};
  \node[dotn] at (-0.180,-0.250) {};
  \node[dotn] at (-0.450,-0.850) {};
  \node[dotn] at (-0.450,0.850) {};
  \node[dotn] at (4.100,0.850) {};
  \node[dotn] at (4.100,-0.850) {};
  \node[dotn] at (3.830,-0.250) {};
  \node[dotn] at (3.470,-0.250) {};
  \node[dotn] at (3.200,-0.850) {};
  \node[dotn] at (3.200,0.850) {};
  \node[dotn] at (7.750,0.850) {};
  \node[dotn] at (7.750,-0.850) {};
  \node[dotn] at (7.480,-0.250) {};
  \node[dotn] at (7.120,-0.250) {};
  \node[dotn] at (6.850,-0.850) {};
  \node[dotn] at (6.850,0.850) {};
  \node[dotn] at (11.400,0.850) {};
  \node[dotn] at (11.400,-0.850) {};
  \node[dotn] at (11.130,-0.250) {};
  \node[dotn] at (10.770,-0.250) {};
  \node[dotn] at (10.500,-0.850) {};
  \node[dotn] at (10.500,0.850) {};
  \node[lab] at (-1.100,0.920) {$I_0^{i}$};
  \node[lab] at (-1.100,0.070) {$I_1^{i}$};
  \node[lab] at (-1.100,-0.780) {$I_2^{i}$};
  \node[lab] at (1.100,0.920) {$O_0^{i}$};
  \node[lab] at (1.100,0.070) {$O_1^{i}$};
  \node[lab] at (1.100,-0.780) {$O_2^{i}$};
  \node[lab] at (2.550,0.920) {$I_0^{i+1}$};
  \node[lab] at (2.550,0.070) {$I_1^{i+1}$};
  \node[lab] at (2.550,-0.780) {$I_2^{i+1}$};
  \node[lab] at (4.750,0.920) {$O_0^{i+1}$};
  \node[lab] at (4.750,0.070) {$O_1^{i+1}$};
  \node[lab] at (4.750,-0.780) {$O_2^{i+1}$};
  \node[lab] at (6.200,0.920) {$I_0^{i+2}$};
  \node[lab] at (6.200,0.070) {$I_1^{i+2}$};
  \node[lab] at (6.200,-0.780) {$I_2^{i+2}$};
  \node[lab] at (8.400,0.920) {$O_0^{i+2}$};
  \node[lab] at (8.400,0.070) {$O_1^{i+2}$};
  \node[lab] at (8.400,-0.780) {$O_2^{i+2}$};
  \node[lab] at (9.850,0.920) {$I_0^{i+3}$};
  \node[lab] at (9.850,0.070) {$I_1^{i+3}$};
  \node[lab] at (9.850,-0.780) {$I_2^{i+3}$};
  \node[lab] at (12.050,0.920) {$O_0^{i+3}$};
  \node[lab] at (12.050,0.070) {$O_1^{i+3}$};
  \node[lab] at (12.050,-0.780) {$O_2^{i+3}$};
\end{tikzpicture}%
}

\newcommand{\rwcycleblock}[2]{%
  \node[draw=black!65,fill=black!3,rounded corners=1.5pt,minimum width=1.30cm,
    minimum height=1.70cm,inner sep=0pt] at (#1,0) {};%
  \node[font=\small] at (#1,0) {$H^{#2}$};%
  \draw[cycsolid] (#1-0.65,0.85) -- (#1-0.85,0.85);%
  \draw[cycsolid] (#1-0.65,-0.85) -- (#1-0.85,-0.85);%
  \draw[cycsolid] (#1+0.65,0.85) -- (#1+0.85,0.85);%
  \draw[cycsolid] (#1+0.65,-0.85) -- (#1+0.85,-0.85);%
  \draw[cycpath] (#1-0.65,0.85) .. controls (#1-0.50,0.60) and (#1-0.50,-0.60) .. (#1-0.65,-0.85);%
  \draw[cycpath] (#1+0.65,0.85) .. controls (#1+0.50,0.60) and (#1+0.50,-0.60) .. (#1+0.65,-0.85);%
  \rwportlabels{#1}{#2}%
}

\newcommand{\ringcyclepanel}{%
\begin{tikzpicture}[x=1cm,y=1cm,line cap=round,line join=round]
  \tikzset{cycsolid/.style={draw=blue!72!black,line width=0.65pt},
    cycunused/.style={draw=gray!70,line width=0.65pt},
    cycpath/.style={draw=blue!72!black,line width=0.8pt,
      dash pattern=on 2.2pt off 1.35pt,line join=bevel},
    cycarrow/.style={->,draw=blue!72!black,line width=0.65pt}}
  \rwlongwires{cycsolid}
  \foreach \x in {0,3.65,7.30}{%
    \draw[cycsolid] (\x+0.85,0.85) -- (\x+2.80,0.85);
  }
  \draw[cycsolid] (-1.40,0.85) -- (-0.85,0.85);
  \draw[cycsolid] (11.80,0.85) -- (13.20,0.85);
  \foreach \x in {0,3.65,7.30}{%
    \draw[cycunused] (\x+0.65,0) -- (\x+3.00,0);
  }
  \draw[cycunused] (-1.40,0) -- (-0.65,0);
  \draw[cycunused] (11.60,0) -- (13.20,0);
  \draw[cycarrow] (1.915,0.85) -- (1.735,0.85);
  \draw[cycarrow] (5.565,0.85) -- (5.385,0.85);
  \draw[cycarrow] (9.210,0.85) -- (9.030,0.85);
  \draw[cycarrow] (12.865,0.85) -- (12.685,0.85);
  \foreach \x in {0,3.65,7.30,10.95}{%
    \draw[cycarrow] (\x-0.09,-1.25) -- (\x+0.09,-1.25);
  }
  \node[font=\small,anchor=west] at (-1.40,1.75) {(d)};
  \rwcycleblock{0}{0}
  \rwcycleblock{3.65}{1}
  \rwcycleblock{7.30}{2}
  \rwcycleblock{10.95}{3}
\end{tikzpicture}%
}

\newcommand{\ringwiring}[1]{%
\ifnum#1=0 \ringopaquepanel\else\ifnum#1=2 \ringcyclepanel\else\ringhexpanel\fi\fi
}

\usepackage[font=small]{caption}
\usepackage[hidelinks]{hyperref}
\hypersetup{pdftitle={Cubic graphs with circumference ratio 21/22 and cyclic connectivity 5},pdfauthor={J\'an Maz\'ak}}
\newtheorem{theorem}{Theorem}
\newtheorem{lemma}[theorem]{Lemma}
\newcommand{\doi}[1]{\url{https://doi.org/#1}}

\title{Cubic graphs with circumference ratio $21/22$\\ and cyclic connectivity $5$}
\author{J\'an Maz\'ak\\[4pt]
Comenius University\\
Mlynsk\'a dolina\\
842 48 Bratislava\\
Slovakia\\[4pt]\texttt{jan.mazak@fmph.uniba.sk}}
\date{}

\begin{document}
\maketitle

\begin{abstract}\sloppy
We construct a cyclically $5$-connected cubic graph of
order $22m$ and circumference $21m$ for every $m\geq 3$. Consequently, the shortness coefficient of the class
of cyclically $5$-connected cubic graphs is at most $21/22$, improving the
previously known upper bound of $24/25$. We also demonstrate how I-reductions,
which delete an edge and suppress its two end-vertices, can be used to prove
cyclic connectivity.
\end{abstract}



\section{The new construction}

For a cubic graph $G$, its \emph{circumference} $\operatorname{circ}(G)$ is the length of a longest
cycle, and its \emph{circumference ratio} is $\operatorname{circ}(G)/|V(G)|$.
For an infinite class of graphs, the \emph{shortness coefficient}
is the lim inf of these ratios as the order tends to infinity.
A cubic graph is cyclically $k$-connected if there is no edge cut of
fewer than $k$ edges both of whose sides contain a cycle.
How long a cycle is forced by high connectivity in a cubic graph is a long-standing open
question. Every $3$-connected cubic graph has a cycle of polynomial length in its order,
although no linear lower bound is known to hold in general~\cite{BJMY,BS,J,LYZ}.
Bondy asked whether every cyclically $4$-connected cubic graph
has a cycle of linear length~\cite[Conjecture~1]{J}, and the dominating cycle
conjecture~\cite[Conjecture~2]{J} would give a
cycle through at least $3/4$ of its vertices. No linear lower bound is known for general
cubic graphs of any fixed cyclic connectivity. We refer to~\cite{LSVZ,MM} for further background.

M\'a\v{c}ajov\'a and Maz\'ak \cite{MM} constructed graphs with circumference ratio
bounded from above by a constant smaller than $1$ for cyclic connectivity in $\{4,5,6\}$.
This paper is only concerned with cyclically $5$-connected cubic graphs,
for which their best ratio was $24/25$.
Their construction is based on repeating a $25$-vertex block that does not admit a 3-edge-colouring.
By Lemma~1 of \cite{MM}, a cycle not contained in that block misses at least one of its vertices.
This gives a local upper bound on circumference. Our construction uses a
smaller $22$-vertex block $H$ and a global argument. Colouring still
plays a central role in the proof, though.

We use the multipole terminology of~\cite[Section~2]{MRS}.
A multipole may contain dangling edges, each with one end
incident with a vertex and the other end free. All multipoles used here are
cubic. A $k$-pole has $k$ dangling edges, and a
$(k_1,\ldots,k_\ell)$-pole has its dangling edges partitioned into connectors
of sizes $k_1,\ldots,k_\ell$. The junction of two dangling edges glues their
free ends into a new edge joining their end vertices, and the junction of
two connectors of the same size performs the junctions of their dangling
edges in a given pairing; we call the edges created in this way junction
edges. A 3-edge-colouring is a proper colouring of all edges, including the dangling
edges.

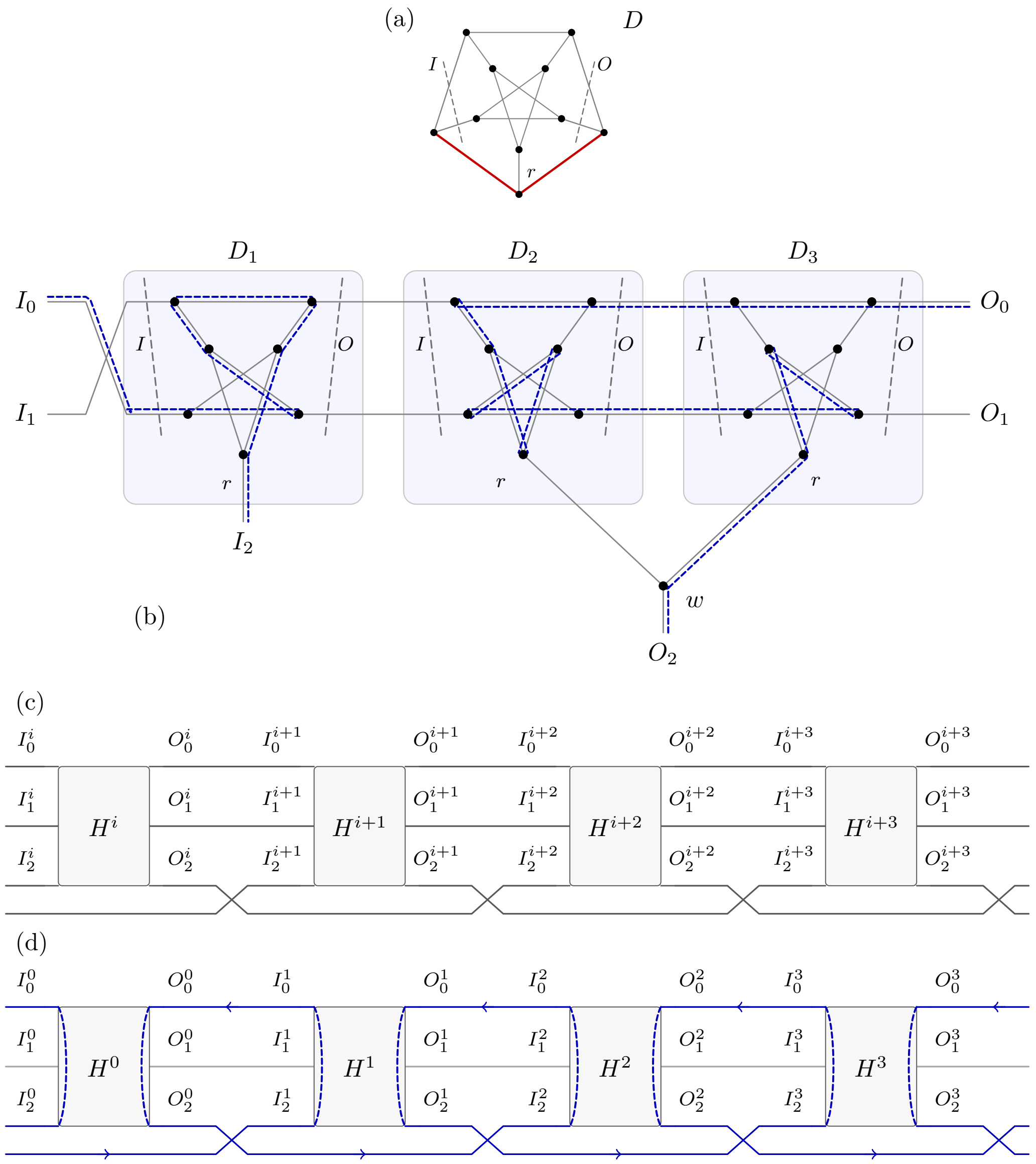
\begin{figure}[p]
\centering
\resizebox{!}{2.6cm}{
\begin{tikzpicture}[x=0.82cm,y=0.82cm,line cap=round,line join=round]
  \tikzset{
    graph edge/.style={draw=black!48,line width=0.55pt},
    restored/.style={draw=red!78!black,line width=1.15pt},
    connector/.style={draw=black!48,line width=0.55pt},
    pvertex/.style={circle,fill=black,draw=black,inner sep=1.15pt,
      minimum size=2.3pt},
    deleted/.style={circle,fill=black,draw=black,inner sep=1.15pt,
      minimum size=2.3pt},
    pair mark/.style={draw=black!55,densely dashed,line width=0.65pt},
    group label/.style={font=\small,inner sep=1pt}
  }

  \node[deleted] (v0) at (-1.902,-0.618) {};
  \node[deleted] (v1) at (0,-2) {};
  \node[deleted] (v2) at (1.902,-0.618) {};
  \node[pvertex] (v3) at (1.176,1.618) {};
  \node[pvertex] (v4) at (-1.176,1.618) {};
  \node[pvertex] (v5) at (-0.951,-0.309) {};
  \node[pvertex] (v6) at (0,-1) {};
  \node[pvertex] (v7) at (0.951,-0.309) {};
  \node[pvertex] (v8) at (0.588,0.809) {};
  \node[pvertex] (v9) at (-0.588,0.809) {};

  \foreach \u/\v in {0/1,1/2,2/3,3/4,4/0,5/7,7/9,9/6,6/8,8/5,0/5,1/6,2/7,3/8,4/9}
    \draw[graph edge] (v\u) -- (v\v);

  \draw[connector] (v0) -- (v4);
  \draw[connector] (v0) -- (v5);
  \draw[connector] (v2) -- (v3);
  \draw[connector] (v2) -- (v7);
  \draw[connector] (v1) -- (v6);

  \draw[pair mark] (-1.69,0.96) -- (-1.26,-0.88);
  \draw[pair mark] (1.69,0.96) -- (1.26,-0.88);

  \draw[restored] (v0) -- (v1) -- (v2);
  \node[group label] at (-1.92,0.92) {$I$};
  \node[group label] at (1.92,0.92) {$O$};
  \node[group label] at (0.28,-1.52) {$r$};
  \node[font=\small,anchor=east] at (-2.15,1.90) {\scalebox{1.29}{(a)}};
  \node[font=\small,anchor=west] at (2.15,1.90) {\scalebox{1.29}{$D$}};

  \foreach \v in {0,1,2}{\node[deleted] at (v\v) {};}
  \foreach \v in {3,4,5,6,7,8,9}{\node[pvertex] at (v\v) {};}
\end{tikzpicture}
}
\par\vspace{14pt}
\resizebox{\textwidth}{!}{
\begin{tikzpicture}[x=0.82cm,y=0.82cm,line cap=round,line join=round]
  \tikzset{
    graph edge/.style={draw=black!48,line width=0.55pt},
    passage/.style={draw=blue!72!black,line width=0.8pt,
      dash pattern=on 2.2pt off 1.35pt,line join=bevel},
    danglinglabel/.style={font=\small,inner sep=4pt},
    block label/.style={font=\small,inner sep=1.5pt},
    vertex dot/.style={circle,fill=black,draw=black,inner sep=1.15pt,
      minimum size=2.3pt},
    pair mark/.style={draw=black!55,densely dashed,line width=0.65pt}
  }

  \filldraw[rounded corners=5pt,fill=blue!4,draw=black!22,line width=0.45pt]
    (-2.05,-1.85) rectangle (2.05,2.15);
  \filldraw[rounded corners=5pt,fill=blue!4,draw=black!22,line width=0.45pt]
    (2.75,-1.85) rectangle (6.85,2.15);
  \filldraw[rounded corners=5pt,fill=blue!4,draw=black!22,line width=0.45pt]
    (7.55,-1.85) rectangle (11.65,2.15);

  \node[block label] at (0,2.47) {$D_1$};
  \node[block label] at (4.8,2.47) {$D_2$};
  \node[block label] at (9.6,2.47) {$D_3$};
  \node[font=\scriptsize,inner sep=0.6pt] at (-1.76,0.90) {$I$};
  \node[font=\scriptsize,inner sep=0.6pt] at (1.76,0.90) {$O$};
  \node[font=\scriptsize,inner sep=0.6pt] at (3.04,0.90) {$I$};
  \node[font=\scriptsize,inner sep=0.6pt] at (6.56,0.90) {$O$};
  \node[font=\scriptsize,inner sep=0.6pt] at (7.84,0.90) {$I$};
  \node[font=\scriptsize,inner sep=0.6pt] at (11.36,0.90) {$O$};

  \coordinate (h0) at (-0.588,0.809);
  \coordinate (h1) at (0.588,0.809);
  \coordinate (h2) at (1.176,1.618);
  \coordinate (h3) at (-0.951,-0.309);
  \coordinate (h4) at (0.951,-0.309);
  \coordinate (h5) at (-1.176,1.618);
  \coordinate (h6) at (0,-1);
  \coordinate (h10) at (5.751,-0.309);
  \coordinate (h11) at (5.976,1.618);
  \coordinate (h12) at (5.388,0.809);
  \coordinate (h13) at (4.212,0.809);
  \coordinate (h14) at (4.8,-1);
  \coordinate (h19) at (3.624,1.618);
  \coordinate (h20) at (3.849,-0.309);
  \coordinate (h7) at (9.6,-1);
  \coordinate (h8) at (8.424,1.618);
  \coordinate (h9) at (8.649,-0.309);
  \coordinate (h15) at (10.188,0.809);
  \coordinate (h16) at (9.012,0.809);
  \coordinate (h17) at (10.551,-0.309);
  \coordinate (h18) at (10.776,1.618);
  \coordinate (h21) at (7.2,-3.25);

  \foreach \u/\v in {
    0/4,0/5,0/6,1/2,1/3,1/6,2/5,2/19,3/4,4/20,
    7/15,7/16,7/21,8/11,8/16,8/18,9/10,9/15,9/17,
    10/13,10/20,11/12,11/19,12/14,12/20,13/14,13/19,14/21,
    15/18,16/17}{\draw[graph edge] (h\u) -- (h\v);}

  \draw[graph edge] (h3) -- (-2.00,-0.309) -- (-2.70,1.618) -- (-3.35,1.618)
    node[danglinglabel,left] {$I_0$};
  \draw[graph edge] (h5) -- (-2.00,1.618) -- (-2.70,-0.309) -- (-3.35,-0.309)
    node[danglinglabel,left] {$I_1$};
  \draw[graph edge] (h6) -- (0,-2.15)
    node[danglinglabel,below] {$I_2$};
  \draw[graph edge] (h18) -- (12.45,1.618)
    node[danglinglabel,right] {$O_0$};
  \draw[graph edge] (h17) -- (12.45,-0.309)
    node[danglinglabel,right] {$O_1$};
  \draw[graph edge] (h21) -- (7.2,-4.05)
    node[danglinglabel,below] {$O_2$};

  \draw[pair mark] (-1.70,2.02) -- (-1.40,-0.72);
  \draw[pair mark] (1.70,2.02) -- (1.40,-0.72);
  \draw[pair mark] (3.10,2.02) -- (3.40,-0.72);
  \draw[pair mark] (6.50,2.02) -- (6.20,-0.72);
  \draw[pair mark] (7.90,2.02) -- (8.20,-0.72);
  \draw[pair mark] (11.30,2.02) -- (11.00,-0.72);

  \foreach \u/\v in {
    3/4,4/0,0/5,5/2,2/1,1/6,
    18/8,8/11,11/19,19/13,13/14,14/12,12/20,20/10,10/9,
    9/17,17/16,16/7,7/21}{\path (h\u) -- (h\v);}

  \draw[passage]
    ([xshift=0pt,yshift=2pt] -3.35,1.618) --
    ([xshift=0pt,yshift=2pt] -2.70,1.618) --
    ([xshift=1.880pt,yshift=0.683pt] -2.70,1.618) --
    ([xshift=1.880pt,yshift=0.683pt] -2.00,-0.309) --
    ([xshift=0pt,yshift=2pt] -2.00,-0.309) --
    ([xshift=0pt,yshift=2pt]h3) --
    ([xshift=0pt,yshift=2pt]h3) --
    ([xshift=0pt,yshift=2pt]h4) --
    ([xshift=-1.18pt,yshift=-1.62pt]h4) --
    ([xshift=-1.18pt,yshift=-1.62pt]h0) --
    ([xshift=-1.62pt,yshift=-1.18pt]h0) --
    ([xshift=-1.62pt,yshift=-1.18pt]h5) --
    ([xshift=0pt,yshift=2pt]h5) --
    ([xshift=0pt,yshift=2pt]h2) --
    ([xshift=1.62pt,yshift=-1.18pt]h2) --
    ([xshift=1.62pt,yshift=-1.18pt]h1) --
    ([xshift=1.90pt,yshift=-0.62pt]h1) --
    ([xshift=1.90pt,yshift=-0.62pt]h6) --
    ([xshift=2pt,yshift=0pt]h6) --
    ([xshift=2pt,yshift=0pt] 0,-2.15);
  \draw[passage]
    ([xshift=0pt,yshift=-2pt] 12.45,1.618) --
    ([xshift=0pt,yshift=-2pt]h18) --
    ([xshift=0pt,yshift=-2pt]h18) --
    ([xshift=0pt,yshift=-2pt]h8) --
    ([xshift=0pt,yshift=-2pt]h8) --
    ([xshift=0pt,yshift=-2pt]h11) --
    ([xshift=0pt,yshift=-2pt]h11) --
    ([xshift=0pt,yshift=-2pt]h19) --
    ([xshift=1.62pt,yshift=1.18pt]h19) --
    ([xshift=1.62pt,yshift=1.18pt]h13) --
    ([xshift=1.90pt,yshift=0.62pt]h13) --
    ([xshift=1.90pt,yshift=0.62pt]h14) --
    ([xshift=-1.90pt,yshift=0.62pt]h14) --
    ([xshift=-1.90pt,yshift=0.62pt]h12) --
    ([xshift=1.18pt,yshift=-1.62pt]h12) --
    ([xshift=1.18pt,yshift=-1.62pt]h20) --
    ([xshift=0pt,yshift=2pt]h20) --
    ([xshift=0pt,yshift=2pt]h10) --
    ([xshift=0pt,yshift=2pt]h10) --
    ([xshift=0pt,yshift=2pt]h9) --
    ([xshift=0pt,yshift=2pt]h9) --
    ([xshift=0pt,yshift=2pt]h17) --
    ([xshift=-1.18pt,yshift=-1.62pt]h17) --
    ([xshift=-1.18pt,yshift=-1.62pt]h16) --
    ([xshift=1.90pt,yshift=0.62pt]h16) --
    ([xshift=1.90pt,yshift=0.62pt]h7) --
    ([xshift=1.37pt,yshift=-1.46pt]h7) --
    ([xshift=1.37pt,yshift=-1.46pt]h21) --
    ([xshift=2pt,yshift=0pt]h21) --
    ([xshift=2pt,yshift=0pt] 7.2,-4.05);

  \node[font=\scriptsize,inner sep=0.6pt] at (-0.28,-1.52) {$r$};
  \node[font=\scriptsize,inner sep=0.6pt] at (4.42,-1.48) {$r$};
  \node[font=\scriptsize,inner sep=0.6pt] at (9.82,-1.46) {$r$};

  \foreach \v in {0,1,2,3,4,5,6,7,8,9,10,11,12,13,14,15,16,17,18,19,20,21}{
    \node[vertex dot] at (h\v) {};
  }
  \node[font=\small,anchor=west] at (7.42,-3.52) {$w$};
  \node[font=\small,anchor=west] at (-2.05,-3.80) {(b)};
\end{tikzpicture}
}
\par\vspace{2pt}
\resizebox{\textwidth}{!}{\ringwiring{0}}
\par\vspace{2pt}
\resizebox{\textwidth}{!}{\ringwiring{2}}
\caption{(a) The Petersen negator $D$, with the deleted path restored in red.
(b) The 6-pole $H$. The dashed blue paths join $I_0$ to $I_2$ and
$O_0$ to $O_2$, covering all but one vertex. (c) A
fragment of $G_m$ showing the junctions between copies of $H$; the pattern
repeats with indices modulo $m$. (d) The cycle of length $21m$, drawn for $m=4$: the two free ends on each row
are joined, and the arrows trace the cycle once.}
\label{fig:block}
\end{figure}

Figure~\ref{fig:block} shows the
22-vertex $(3,3)$-pole
$H$, with ordered input connector $(I_0,I_1,I_2)$ and output connector
$(O_0,O_1,O_2)$.
The graph6 code of $H$ is
\[
\texttt{\detokenize{UIee??????_C?@?A??WA_?W??OG@?}\textasciigrave{}\detokenize{?I?GI??OG?}},
\]
with dangling edges attached at vertices $3,5,6,18,17,21$.

Structurally, $H$ consists of three copies $D_1,D_2,D_3$ of the Petersen negator $D$,
together with a vertex $w$. The Petersen negator is obtained by deleting the
three vertices on the red path in Figure~\ref{fig:block}(a), leaving a dangling
edge at each surviving neighbour \cite[Sections~4 and~5.1]{MRS}.

The dangling edges at the two deleted end vertices form the 2-connectors
$I$ and $O$; the fifth is the residual dangling edge $r$.
Thus $D=D(I,O,r)$ is a $(2,2,1)$-pole. Since the Petersen graph has girth 5,
the five dangling edges of $D$ are incident with distinct vertices, and the
two vertices incident with each 2-connector are nonadjacent.
Join the three negators and $w$ as in Figure~\ref{fig:block}(b), with the
ordered dangling edges shown there.

For each integer
$m\geq3$, take $m$ copies $H^i$, with indices modulo $m$. In $H^i$, write
$I_j^i,O_j^i$ for the dangling edges corresponding to $I_j,O_j$. Perform the
junctions
\[
 O_0^i I_0^{i+1},\qquad O_1^i I_1^{i+1},\qquad
 O_2^i I_2^{i+2}.
\]
The junction edges of the first two types are \emph{short} and those of the third
type are \emph{long}; accordingly, $I_0^i,I_1^i,O_0^i,O_1^i$ are \emph{short dangling
edges} and $I_2^i,O_2^i$ are \emph{long dangling edges} before the junctions.
Call the resulting graph $G_m$. Clearly, $G_m$ is a simple connected cubic
graph of order $22m$. The restriction $m\geq3$ is needed: the four
short junction edges between $H^0$ and $H^1$ would form a $4$-cut in $G_2$.

\begin{theorem}
\label{thm:main}
For every integer $m\geq3$, the graph $G_m$ is cyclically $5$-connected,
and has $22m$ vertices and circumference
\[
 \operatorname{circ}(G_m)=21m=\tfrac{21}{22}\,|V(G_m)|.
\]
\end{theorem}
Lemmas
\ref{lem:long-edge} and~\ref{lem:no-long-edge} give the circumference upper
bound, while Lemma~\ref{lem:lower-bound} gives the matching lower bound.
Cyclic $5$-connectivity is proved in Section~\ref{sec:cyclic-connectivity} using I-reductions.

\section{Circumference bound}

A \emph{path system} in a multipole is a set of vertex-disjoint paths,
each with at least one edge, whose ends lie on dangling edges. A vertex is
\emph{covered} if it lies on one of the paths. We first record two
colouring arguments that we use repeatedly.

\begin{lemma}
\label{lem:path-colouring}
Let $P$ be a spanning path system in a cubic multipole. Colour each path, including its two end dangling edges, alternately with colours 1 and 2, and all unused edges with colour 3.
This is a proper 3-edge-colouring, and colours 1 and 2 may be interchanged
independently on each path.
\end{lemma}
\begin{proof}
As in~\cite[Lemma~1]{MM}, every vertex has two path edges receiving
colours 1 and 2 and one unused edge receiving colour 3. Interchanging
colours along a path preserves this property.
\end{proof}

Let $Q$ be the $15$-vertex $(2,2,1)$-pole obtained from $H$ by deleting
the vertices of $D_1$, keeping the two dangling edges at $D_2$.

\begin{lemma}
\label{lem:piece-colourings}
In every 3-edge-colouring of $D$, the two edges of exactly one 2-connector have the same colour, and the residual dangling edge has the colour missing from the other 2-connector. In every 3-edge-colouring of $Q$, the two edges of each 2-connector have different colours. In every 3-edge-colouring of $H$, $I_2$ and $O_2$ have
different colours.
\end{lemma}
\begin{proof}
The assertion for $D$ is the negator colouring rule \cite{MS,MRS}.
If the two edges between $D_2$ and $D_3$ in $Q$ had different colours, the residual edges of $D_2$ and $D_3$ would both have the colour missing from them, contradicting properness at $w$. Thus these two edges have the same colour, and the two edges of each outer 2-connector of $Q$ have different colours. In $H$, the two edges between $D_1$ and $D_2$ therefore have different colours. The colour missing from them occurs on the residual edges of $D_1$ and $D_2$,
namely $I_2$ and the edge joining $D_2$ to $w$. Properness at $w$ forces $I_2$ and $O_2$ to have different colours.
\end{proof}

To bound cycles which use a long edge, we work with chains.
A \emph{chain} is a sequence $H^0,\ldots,H^{L-1}$ of copies of $H$ joined
by the short junction edges $O_a^{i-1}I_a^i$ for $a\in\{0,1\}$ and
$1\leq i<L$. Its dangling edges are all the long ones and also the short ones at
its two ends. A single copy of $H$ is a chain with one block.

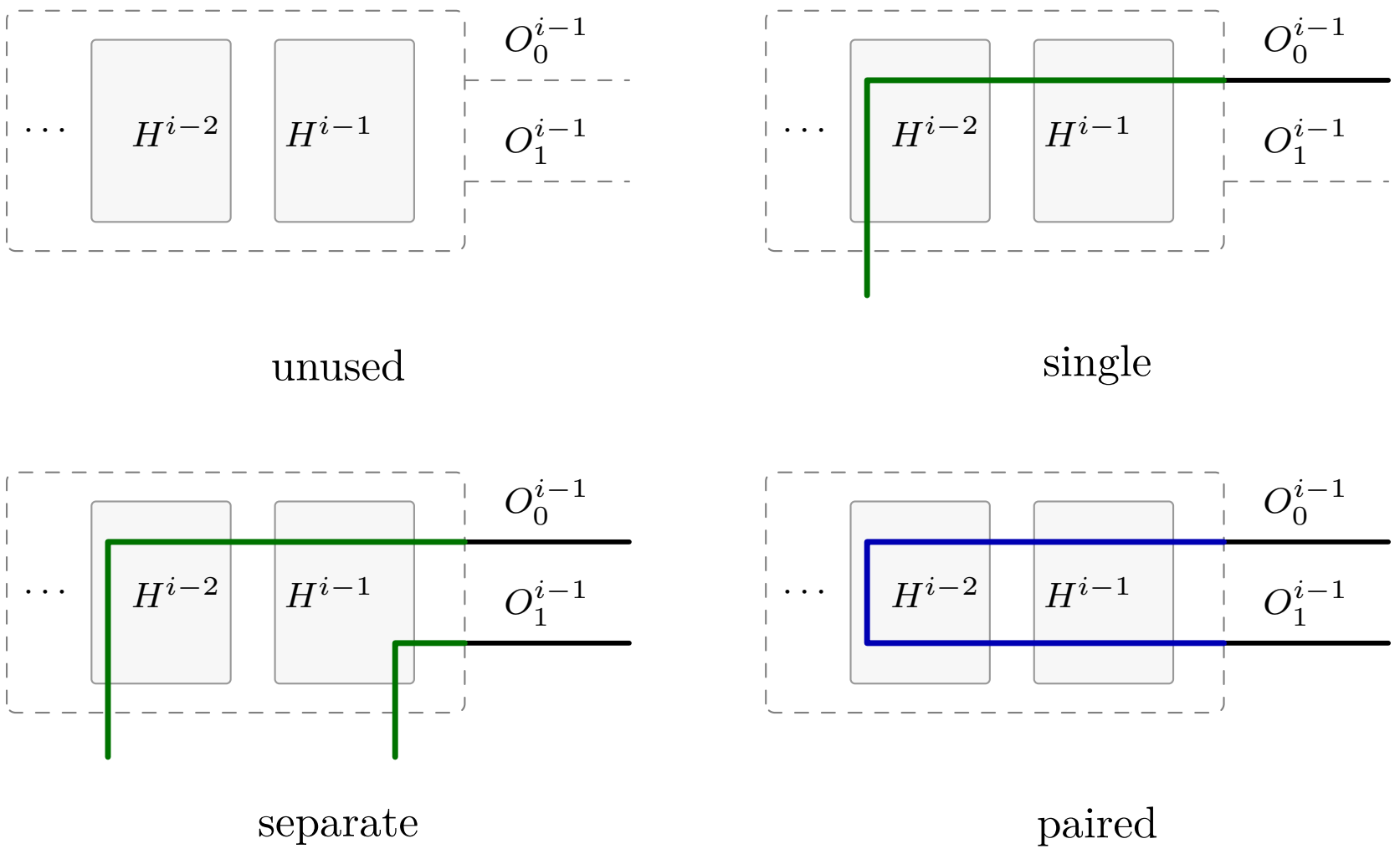
\begin{figure}[H]
\centering
\newcommand{\gsframe}[1]{%
  \draw[leftpart] (-3.62,-0.95) rectangle (0,0.95);
  \draw[blockbox] (-2.95,-0.72) rectangle (-1.85,0.72);
  \draw[blockbox] (-1.50,-0.72) rectangle (-0.40,0.72);
  \node[gslab] at (-2.28,0) {$H^{i-2}$};
  \node[gslab] at (-1.07,0) {$H^{i-1}$};
  \node[gslab] at (-3.29,0) {$\cdots$};
  \node[gslab,anchor=south] at (0.65,0.53) {$O_0^{i-1}$};
  \node[gslab,anchor=south] at (0.65,-0.27) {$O_1^{i-1}$};
  \node[align=center] at (-1.0,-1.85) {#1};
}
\begin{tikzpicture}[x=1cm,y=1cm,line cap=round,line join=round,every node/.style={font=\small}]
  \tikzset{
    cedge/.style={black,line width=1.1pt},
    notc/.style={gray,dashed},
    longpath/.style={green!45!black,line width=1.3pt},
    pairedpath/.style={blue!70!black,line width=1.3pt},
    leftpart/.style={gray,dashed,rounded corners=2pt},
    blockbox/.style={draw=black!40,fill=black!3,rounded corners=1pt},
    gslab/.style={font=\scriptsize,inner sep=0pt}
  }
  \begin{scope}[shift={(0,0)}]
    \gsframe{unused}
    \draw[notc] (0,0.4) -- (1.3,0.4);
    \draw[notc] (0,-0.4) -- (1.3,-0.4);
  \end{scope}
  \begin{scope}[shift={(6.0,0)}]
    \gsframe{single}
    \draw[cedge] (0,0.4) -- (1.3,0.4);
    \draw[notc] (0,-0.4) -- (1.3,-0.4);
    \draw[longpath] (0,0.4) -- (-2.82,0.4) -- (-2.82,-1.30);
  \end{scope}
  \begin{scope}[shift={(0,-3.65)}]
    \gsframe{separate}
    \draw[cedge] (0,0.4) -- (1.3,0.4);
    \draw[cedge] (0,-0.4) -- (1.3,-0.4);
    \draw[longpath] (0,0.4) -- (-2.82,0.4) -- (-2.82,-1.30);
    \draw[longpath] (0,-0.4) -- (-0.55,-0.4) -- (-0.55,-1.30);
  \end{scope}
  \begin{scope}[shift={(6.0,-3.65)}]
    \gsframe{paired}
    \draw[cedge] (0,0.4) -- (1.3,0.4);
    \draw[cedge] (0,-0.4) -- (1.3,-0.4);
    \draw[pairedpath] (0,0.4) -- (-2.82,0.4) -- (-2.82,-0.4) -- (0,-0.4);
  \end{scope}
\end{tikzpicture}

\caption{The four states of $J_i$. Solid lines are edges of $P$. The dashed
box represents the blocks $H^0,\ldots,H^{i-1}$, and long edges point downwards.}
\label{fig:gap-states}
\end{figure}

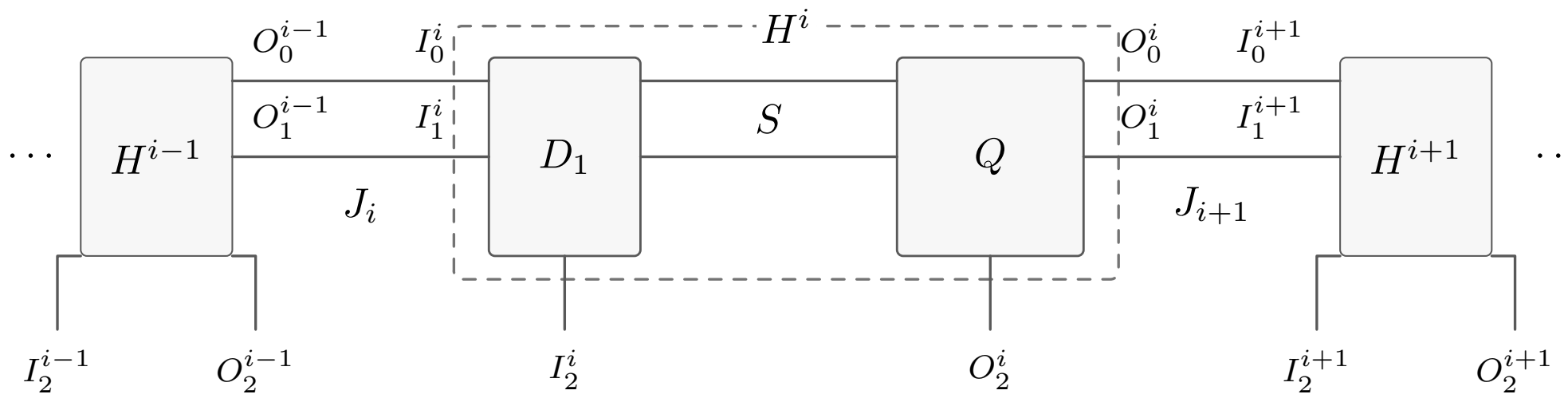
\begin{figure}[H]
\centering
\resizebox{\textwidth}{!}{
\newcommand{\rgbox}[2]{%
  \node[draw=black!65,fill=black!3,rounded corners=1.5pt,minimum width=1.30cm,
    minimum height=1.70cm,inner sep=0pt] at (#1,0) {};%
  \node[font=\small,inner sep=0pt] at (#1,0) {#2};%
}
\newcommand{\rglong}[3]{%
  \draw[rgwire] (#1,-0.85) -- (#2,-0.85) -- (#2,-1.48);%
  \node[font=\scriptsize,anchor=north,inner sep=0pt] at (#2,-1.65) {#3};%
}
\begin{tikzpicture}[x=1cm,y=1cm,line cap=round,line join=round]
  \tikzset{
    rgwire/.style={draw=black!65,line width=0.65pt},
    rglab/.style={font=\scriptsize,anchor=south,inner sep=0pt}
  }
  \rgbox{0}{$H^{i-1}$}
  \rgbox{10.8}{$H^{i+1}$}

  \draw[draw=black!55,dashed,rounded corners=2pt,line width=0.65pt]
    (2.55,-1.05) rectangle (8.25,1.12);
  \node[font=\small,inner sep=1pt,fill=white] at (5.40,1.12) {$H^i$};
  \draw[draw=black!65,fill=black!3,rounded corners=1.5pt,line width=0.55pt]
    (2.85,-0.85) rectangle (4.15,0.85);
  \node[font=\small,inner sep=0pt] at (3.50,0) {$D_1$};
  \draw[draw=black!65,fill=black!3,rounded corners=1.5pt,line width=0.55pt]
    (6.35,-0.85) rectangle (7.95,0.85);
  \node[font=\small,inner sep=0pt] at (7.15,0) {$Q$};

  \draw[rgwire] (0.65,0.65) -- (2.85,0.65);
  \draw[rgwire] (0.65,0) -- (2.85,0);
  \node[rglab] at (1.16,0.80) {$O_0^{i-1}$};
  \node[rglab] at (1.16,0.15) {$O_1^{i-1}$};
  \node[rglab] at (2.35,0.80) {$I_0^{i}$};
  \node[rglab] at (2.35,0.15) {$I_1^{i}$};
  \node[font=\small,inner sep=0pt] at (1.75,-0.42) {$J_i$};

  \draw[rgwire] (4.15,0.65) -- (6.35,0.65);
  \draw[rgwire] (4.15,0) -- (6.35,0);
  \node[font=\small,inner sep=1pt,fill=white] at (5.25,0.325) {$S$};

  \draw[rgwire] (7.95,0.65) -- (10.15,0.65);
  \draw[rgwire] (7.95,0) -- (10.15,0);
  \node[rglab] at (8.45,0.80) {$O_0^{i}$};
  \node[rglab] at (8.45,0.15) {$O_1^{i}$};
  \node[rglab] at (9.55,0.80) {$I_0^{i+1}$};
  \node[rglab] at (9.55,0.15) {$I_1^{i+1}$};
  \node[font=\small,inner sep=0pt] at (9.05,-0.42) {$J_{i+1}$};

  \rglong{-0.65}{-0.85}{$I_2^{i-1}$}
  \rglong{0.65}{0.85}{$O_2^{i-1}$}
  \rglong{3.50}{3.50}{$I_2^{i}$}
  \rglong{7.15}{7.15}{$O_2^{i}$}
  \rglong{10.15}{9.95}{$I_2^{i+1}$}
  \rglong{11.45}{11.65}{$O_2^{i+1}$}
  \node[font=\small,inner sep=0pt] at (-1.05,0) {$\cdots$};
  \node[font=\small,inner sep=0pt] at (12.05,0) {$\cdots$};
\end{tikzpicture}
}
\caption{Blocks $H^{i-1},H^i,H^{i+1}$ of a chain, with $H^i$ split into
$D_1$ and $Q$ joined by the edge pair $S$.}
\label{fig:ring}
\end{figure}

Let $P$ be a path system in a chain. An edge is \emph{used by $P$} if it
lies on one of its paths. For $1\leq i<L$, let $J_i$ be the
pair of short junction edges joining $H^{i-1}$ to $H^i$; let $J_0$ be the
pair $I_0^0,I_1^0$ and $J_L$ the pair $O_0^{L-1},O_1^{L-1}$.
The \emph{state} of $J_i$ is one of the following (see
Figures~\ref{fig:gap-states} and~\ref{fig:ring}):
{\setlength{\topsep}{2pt}\setlength{\partopsep}{0pt}
\begin{itemize}\setlength{\itemsep}{0pt}\setlength{\parskip}{0pt}
\item \emph{unused} if $P$ contains neither edge of $J_i$;
\item \emph{single} if $P$ contains exactly one edge of $J_i$;
\item \emph{paired} if $P$ contains both edges and the part of $P$ in
$H^0,\ldots,H^{i-1}$ has a path joining them;
\item \emph{separate} if $P$ contains both edges and $J_i$ is not paired.
\end{itemize}}

Fix a block $H^i$ and let $S$ be the pair of edges between $D_1$ and $Q$.
Define the state of $S$ by the same four cases as for $J_i$, with the earlier
part now consisting of $H^0,\ldots,H^{i-1}$ together with $D_1$.

\begin{lemma}
\label{lem:full-negator}
If $P$ covers all seven vertices of $D_1$ and $J_i$ is single or paired,
then $S$ is unused or paired.
\end{lemma}
\begin{proof}
Apply Lemma~\ref{lem:path-colouring} to the paths in $D_1$.
If $J_i$ is single, the two input edges have different colours, so the two edges of $S$ have the same colour by Lemma~\ref{lem:piece-colourings}. Thus $S$ is unused, or its two used edges lie on one path of $D_1$; otherwise interchanging the colours on one path would give them different colours. Hence $S$ is unused or paired.

If $J_i$ is paired, its two edges lie on distinct paths in $D_1$,
since a path joining them would close a cycle with the earlier path.
By Lemma~\ref{lem:path-colouring}, we can interchange colours 1 and 2 on one
of these two paths so that the two edges of $J_i$ receive different colours;
Lemma~\ref{lem:piece-colourings} then forces the two edges of $S$ to have the
same colour. Moreover, the two paths use at least four of the five dangling
edges of $D_1$, so at least one edge of $S$ is used. Since any unused edge has
colour 3 while used edges have colour 1 or 2, both edges of $S$ must be used.
As the number of used dangling edges is even, $I_2$ is unused, so there are
only these two paths in $D_1$, and they connect the edges of $J_i$ to distinct
edges of $S$. Together with the earlier path between the two edges of $J_i$, they show that
$S$ is paired.
\end{proof}

\begin{lemma}
\label{lem:paths-in-Q}
If $P$ covers all 15 vertices of $Q$, then $S$ is single or separate and
$J_{i+1}$ is single or paired.
\end{lemma}
\begin{proof}
Apply Lemma~\ref{lem:path-colouring} to the paths in $Q$. Neither
2-connector can be unused, since by Lemma~\ref{lem:piece-colourings} the
two edges of each have different colours. Let $T$ be either 2-connector
of $Q$. If both edges of $T$ are used but lie on different paths in $Q$,
interchange colours 1 and 2 on one of these paths. By
Lemma~\ref{lem:path-colouring}, the colouring remains proper, but the two
edges of $T$ now have the same colour, contrary to
Lemma~\ref{lem:piece-colourings}. Thus the edges of $T$ are joined by a
path inside $Q$ whenever both are used. For $T=J_{i+1}$, since it is not
unused, this means $J_{i+1}$ is single or paired. For $T=S$, if $S$ were
paired, there would also be a path joining its two edges through the earlier
part, and these two paths would form a cycle in $P$; since $S$ is not unused,
it is single or separate.
\end{proof}

\begin{lemma}
\label{lem:nonempty-Q}
If $S$ is unused or paired and $P$ covers at least one vertex of $Q$,
then $J_{i+1}$ is single or paired.
\end{lemma}
\begin{proof}
Consider the restriction of $P$ to $H^0,\ldots,H^i$. If $S$ is unused, any
path of this restriction that meets $Q$ stays inside $Q$. If $S$ is paired,
the two edges of $S$ are joined by a path in $H^0,\ldots,H^{i-1},D_1$, so any
path meeting $Q$ that exits through one edge of $S$ returns to $Q$ through the
other. In either case, every path of the restriction that meets $Q$ has both
of its ends among $O_0^i,O_1^i,O_2^i$. Since $P$ covers at least one vertex of
$Q$ and only three ends are available, there is a unique such path. It joins
$O_0^i$ to $O_1^i$ (so $J_{i+1}$ is paired) or one of them to $O_2^i$ (so
$J_{i+1}$ is single).
\end{proof}

\begin{lemma}
\label{lem:block}
Let $P$ be a path system in a chain. Then between any two blocks in which
$P$ covers all $22$ vertices there is a block in which $P$ covers at most
$20$ vertices.
\end{lemma}
\begin{proof}
We use two local properties, which we derive from the previous lemmas:
\begin{itemize}
\item[(F1)] If $P$ covers all $22$ vertices of $H^i$, then $J_i$ is unused
or separate and $J_{i+1}$ is single or paired.
\item[(F2)] If $P$ covers exactly $21$ vertices of $H^i$ and $J_i$ is
single or paired, then $J_{i+1}$ is single or paired.
\end{itemize}
If $P$ covers all 22 vertices of $H^i$, Lemma~\ref{lem:paths-in-Q} gives
$S$ single or separate and $J_{i+1}$ single or paired.
Lemma~\ref{lem:full-negator} then excludes a single or paired $J_i$.
This proves (F1). For (F2), if $P$ covers all of $Q$, use
Lemma~\ref{lem:paths-in-Q}. Otherwise its one uncovered vertex lies in
$Q$, so it covers all of $D_1$ and 14 vertices of $Q$.
Lemmas~\ref{lem:full-negator} and~\ref{lem:nonempty-Q} give the required
single or paired output.

Consider two consecutive fully covered blocks. If every block between
them covered 21 vertices, (F1) and (F2) would keep the intervening
junctions single or paired, contrary to (F1) at the second full block.
Thus some block between them has at most 20 vertices covered.
\end{proof}

\begin{lemma}
\label{lem:long-edge}
Every cycle of $G_m$ that uses a long edge has length at most $21m$.
\end{lemma}
\begin{proof}
Let $C$ be such a cycle. Cut every long edge of $G_m$ and the two short edges $O_a^{m-1}I_a^0$, for $a\in\{0,1\}$, keeping both halves of each cut edge as dangling edges. This gives a chain of $m$ blocks, and, since $C$ uses a long edge, the same cuts turn $C$ into a path system in it. Take two disjoint copies of this chain and its path system. For each $a$, join $O_a^{m-1}$ of the first copy to $I_a^0$ of the second, and include the new edge in the path system exactly when $O_a^{m-1}I_a^0$ lies in $C$. Let $P$ be the subgraph formed by the two copies of the path system and the new edges included. It lies in the chain of $2m$ blocks obtained by joining the two copies, and it has exactly $2|V(C)|$ vertices.

The subgraph $P$ contains no cycle. Such a cycle would use both edges joining the two copies, and removing them would leave a path joining $O_0^{m-1}$ and $O_1^{m-1}$ in the first copy and a path joining $I_0^0$ and $I_1^0$ in the second. These correspond to two distinct paths of the path system obtained from $C$ by the cuts, so restoring the two short edges gives a cycle $C$ with no long edge, contradicting the choice of $C$. Thus the components of $P$ form a path system. By Lemma~\ref{lem:block}, the blocks in which $P$ covers all $22$ vertices outnumber those in which it covers at most $20$ by at most one, so
\[
2|V(C)|=|V(P)|=\sum_{i=0}^{2m-1}|V(P)\cap V(H^i)|\leq 42m+1.
\]
Since $|V(C)|$ is an integer, $|V(C)|\leq 21m$.
\end{proof}

We now deal with cycles that use no long edge, using the colouring idea of~\cite{MM}.

\begin{lemma}
\label{lem:no-long-edge}
Every cycle of $G_m$ that uses no long edge has length at most $21m$.
\end{lemma}
\begin{proof}
A cycle contained in one copy of $H$ has length at most $22$, less than $21m$.
Otherwise, the restriction of the cycle to each copy of $H$ is a path system with all ends on $I_0,I_1,O_0,O_1$. Such a path system cannot cover all 22 vertices: by Lemma~\ref{lem:path-colouring} it would give a 3-edge-colouring of $H$ in which the unused dangling edges $I_2$ and $O_2$ both have colour 3, contrary to Lemma~\ref{lem:piece-colourings}. So the cycle misses a vertex in every copy, and it has length at most $21m$.
\end{proof}

\begin{lemma}
\label{lem:lower-bound}
For every $m\geq3$, the graph $G_m$ contains a cycle of length $21m$.
\end{lemma}
\begin{proof}
In every copy take the two paths of Figure~\ref{fig:block}(b), joining $I_0$ to $I_2$ and $O_0$ to $O_2$; they cover 21 vertices, and $I_1,O_1$ are not used by the cycle. Figure~\ref{fig:block}(d) shows the resulting cycle: it uses the short junctions $O_0^iI_0^{i+1}$ and the long junctions $O_2^iI_2^{i+2}$, and traverses both paths of every copy. Starting at $I_0^i$, the walk $I_0^i\to O_0^{i-1}\to O_2^{i-1}\to I_2^{i+1}\to I_0^{i+1}$, which runs through the chosen paths of $H^{i-1}$ and $H^{i+1}$, advances the index by one. So after $m$ rounds it closes up, having used each of the $2m$ chosen paths once; its vertices are the $21m$ vertices of these paths, and it is a cycle of length $21m$.
\end{proof}
\begin{proof}[Proof of Theorem~\ref{thm:main}]
Lemmas~\ref{lem:long-edge} and~\ref{lem:no-long-edge} give $\operatorname{circ}(G_m)\le21m$, and
Lemma~\ref{lem:lower-bound} gives a cycle of length $21m$, so
$\operatorname{circ}(G_m)=21m=\tfrac{21}{22}|V(G_m)|$ for every $m\geq3$. Cyclic $5$-connectivity is proved in Section~\ref{sec:cyclic-connectivity}.
\end{proof}

\section{Cyclic connectivity}
\label{sec:cyclic-connectivity}

In this section we prove that $G_m$ is cyclically $5$-connected for every $m\geq3$.

We chose a picturesque way of proving this: via I-reductions. 
An \emph{I-extension} of a cubic graph
subdivides two distinct edges by new vertices $x$ and $y$ and joins $x$ to $y$
by a \emph{new edge} $xy$. The inverse operation, an \emph{I-reduction},
deletes the edge $xy$ and suppresses its end-vertices $x$ and $y$. The property of I-reduction we are interested in is captured
by the I-extension lemma of M\'a\v{c}ajov\'a and Maz\'ak (Lemma~7 in~\cite{MMpm})
which we restate here for cyclic connectivity $5$. The \emph{distance} between
two edges is the smallest distance between an end of one and an end of the other.

\begin{lemma}[\cite{MMpm}]\label{lem:iext}
Let $G'$ arise by an I-extension from a cyclically $5$-connected cubic graph
$G$ by subdividing two edges at distance at least $2$.
Then $G'$ is cyclically $5$-connected.
\end{lemma}

We call an I-reduction \emph{valid} if its two new edges $e$ and $f$ are
at distance at least $2$ in the reduced graph. In the figures, $e$ and $f$ are the two green dashed edges created by a reduction. Their distance can be read from the drawings together with the junctions of Section~1, also for small $m$, when blocks wrap around.

\begin{proof}[Proof that $G_m$ is cyclically $5$-connected]
All reductions shown in Figure~\ref{fig:iext-block} are valid. Applying them
in every block, panel by panel, takes $G_m$ to the ring $R_m$ of hexagons
in panel (g).

\begin{figure}[!ht]
\centering
\resizebox{\textwidth}{!}{\begin{minipage}{17.2cm}\centering
\begin{tikzpicture}[x=0.82cm,y=0.82cm,line cap=round,line join=round]
\tikzset{
    graph edge/.style={draw=black!48,line width=0.55pt},
    reduced/.style={draw=red!78!black,line width=1.8pt},
    new edge/.style={draw=green!50!black,line width=0.9pt,dash pattern=on 2.2pt off 1.35pt},
    vertex dot/.style={circle,fill=black,draw=black,inner sep=1.0pt,minimum size=2.0pt},
    vlabel/.style={font=\scriptsize,text=black!85,inner sep=0.3pt},
    port/.style={font=\scriptsize,inner sep=0.6pt},
    stage/.style={font=\scriptsize,text=red!78!black,inner sep=0.4pt}
  }
\node[font=\small,anchor=west] at (-2.900,2.600) {(a)};
\draw[graph edge] (-0.588,0.809) -- (0.951,-0.309);
\draw[graph edge] (-0.588,0.809) -- (-1.176,1.618);
\draw[graph edge] (0.588,0.809) -- (1.176,1.618);
\draw[graph edge] (0.588,0.809) -- (-0.951,-0.309);
\draw[graph edge] (0.588,0.809) -- (0.000,-1.000);
\draw[graph edge] (1.176,1.618) -- (-1.176,1.618);
\draw[graph edge] (1.176,1.618) -- (3.624,1.618);
\draw[graph edge] (-0.951,-0.309) -- (0.951,-0.309);
\draw[graph edge] (0.951,-0.309) -- (3.849,-0.309);
\draw[graph edge] (9.600,-1.000) -- (9.012,0.809);
\draw[graph edge] (9.600,-1.000) -- (7.200,-3.250);
\draw[graph edge] (8.424,1.618) -- (5.976,1.618);
\draw[graph edge] (8.424,1.618) -- (9.012,0.809);
\draw[graph edge] (8.424,1.618) -- (10.776,1.618);
\draw[graph edge] (8.649,-0.309) -- (5.751,-0.309);
\draw[graph edge] (8.649,-0.309) -- (10.188,0.809);
\draw[graph edge] (8.649,-0.309) -- (10.551,-0.309);
\draw[graph edge] (5.751,-0.309) -- (4.212,0.809);
\draw[graph edge] (5.976,1.618) -- (5.388,0.809);
\draw[graph edge] (5.976,1.618) -- (3.624,1.618);
\draw[graph edge] (5.388,0.809) -- (4.800,-1.000);
\draw[graph edge] (5.388,0.809) -- (3.849,-0.309);
\draw[graph edge] (4.212,0.809) -- (4.800,-1.000);
\draw[graph edge] (4.212,0.809) -- (3.624,1.618);
\draw[graph edge] (4.800,-1.000) -- (7.200,-3.250);
\draw[graph edge] (10.188,0.809) -- (10.776,1.618);
\draw[graph edge] (9.012,0.809) -- (10.551,-0.309);
\draw[reduced] (9.600,-1.000) -- (10.188,0.809);
\draw[reduced] (-0.588,0.809) -- (0.000,-1.000);
\draw[reduced] (5.751,-0.309) -- (3.849,-0.309);
\draw[graph edge] (-0.951,-0.309) -- (-1.727,-0.309);
\draw[graph edge] (-1.176,1.618) -- (-1.763,1.618);
\draw[graph edge] (0.000,-1.000) -- (0.000,-1.900);
\draw[graph edge] (10.551,-0.309) -- (11.309,-0.309);
\draw[graph edge] (10.776,1.618) -- (11.375,1.618);
\draw[graph edge] (7.200,-3.250) -- (7.200,-4.150);
\node[port,anchor=east] at (-1.910,-0.309) {$I_0$};
\node[port,anchor=east] at (-1.946,1.618) {$I_1$};
\node[port,anchor=north] at (0.000,-2.071) {$I_2$};
\node[port,anchor=west] at (11.492,-0.309) {$O_1$};
\node[port,anchor=west] at (11.558,1.618) {$O_0$};
\node[port,anchor=north] at (7.200,-4.321) {$O_2$};
\node[vertex dot] at (-0.588,0.809) {};
\node[vertex dot] at (0.588,0.809) {};
\node[vertex dot] at (1.176,1.618) {};
\node[vertex dot] at (-0.951,-0.309) {};
\node[vertex dot] at (0.951,-0.309) {};
\node[vertex dot] at (-1.176,1.618) {};
\node[vertex dot] at (0.000,-1.000) {};
\node[vertex dot] at (9.600,-1.000) {};
\node[vertex dot] at (8.424,1.618) {};
\node[vertex dot] at (8.649,-0.309) {};
\node[vertex dot] at (5.751,-0.309) {};
\node[vertex dot] at (5.976,1.618) {};
\node[vertex dot] at (5.388,0.809) {};
\node[vertex dot] at (4.212,0.809) {};
\node[vertex dot] at (4.800,-1.000) {};
\node[vertex dot] at (10.188,0.809) {};
\node[vertex dot] at (9.012,0.809) {};
\node[vertex dot] at (10.551,-0.309) {};
\node[vertex dot] at (10.776,1.618) {};
\node[vertex dot] at (3.624,1.618) {};
\node[vertex dot] at (3.849,-0.309) {};
\node[vertex dot] at (7.200,-3.250) {};
\node[vlabel] at (-0.312,1.021) {0};
\node[vlabel] at (0.312,1.021) {1};
\node[vlabel] at (1.391,1.899) {2};
\node[vlabel] at (-1.223,-0.037) {3};
\node[vlabel] at (1.223,-0.581) {4};
\node[vlabel] at (-0.999,1.924) {5};
\node[vlabel] at (0.276,-1.212) {6};
\node[vlabel] at (9.648,-1.363) {7};
\node[vlabel] at (8.209,1.899) {8};
\node[vlabel] at (8.377,-0.581) {9};
\node[vlabel] at (6.040,-0.598) {10};
\node[vlabel] at (6.217,1.932) {11};
\node[vlabel] at (5.059,1.061) {12};
\node[vlabel] at (4.541,1.061) {13};
\node[vlabel] at (4.743,-1.435) {14};
\node[vlabel] at (9.859,1.061) {15};
\node[vlabel] at (9.341,1.061) {16};
\node[vlabel] at (10.840,-0.020) {17};
\node[vlabel] at (10.587,1.945) {18};
\node[vlabel] at (3.383,1.932) {19};
\node[vlabel] at (3.560,-0.598) {20};
\node[vlabel] at (7.606,-3.359) {21};
\end{tikzpicture}
\par\medskip
\begin{minipage}[t]{0.49\textwidth}\centering
\begin{tikzpicture}[x=0.52cm,y=0.52cm,line cap=round,line join=round]
\tikzset{graph edge/.style={draw=black!48,line width=0.55pt},reduced/.style={draw=red!78!black,line width=1.8pt},new edge/.style={draw=green!50!black,line width=0.9pt,dash pattern=on 2.2pt off 1.35pt},overlap edge/.style={draw=green!50!black,line width=0.65pt,dash pattern=on 2.2pt off 1.35pt},vertex dot/.style={circle,fill=black,draw=black,inner sep=1.0pt,minimum size=2.0pt},vlabel/.style={font=\scriptsize,text=black!85,inner sep=0.3pt},port/.style={font=\scriptsize,inner sep=0.6pt}}
\node[font=\small,anchor=west] at (-2.450,3.050) {(b)};
\draw[graph edge] (0.588,0.809) -- (1.176,1.618);
\draw[graph edge] (0.588,0.809) -- (-0.951,-0.309);
\draw[graph edge] (-0.951,-0.309) -- (0.951,-0.309);
\draw[graph edge] (1.176,1.618) -- (-1.176,1.618);
\draw[graph edge] (8.424,1.618) -- (5.976,1.618);
\draw[graph edge] (5.976,1.618) -- (5.388,0.809);
\draw[graph edge] (5.388,0.809) -- (4.800,-1.000);
\draw[graph edge] (4.212,0.809) -- (4.800,-1.000);
\draw[graph edge] (8.424,1.618) -- (9.012,0.809);
\draw[graph edge] (8.649,-0.309) -- (10.551,-0.309);
\draw[graph edge] (9.012,0.809) -- (10.551,-0.309);
\draw[reduced] (8.424,1.618) -- (10.776,1.618);
\draw[graph edge] (1.176,1.618) -- (3.624,1.618);
\draw[graph edge] (5.976,1.618) -- (3.624,1.618);
\draw[graph edge] (4.212,0.809) -- (3.624,1.618);
\draw[graph edge] (4.800,-1.000) -- (7.200,-3.250);
\draw[graph edge] (-0.951,-0.309) -- (-1.727,-0.309);
\draw[graph edge] (-1.176,1.618) -- (-1.763,1.618);
\draw[graph edge] (10.776,1.618) -- (11.375,1.618);
\draw[graph edge] (10.551,-0.309) -- (11.309,-0.309);
\draw[graph edge] (7.200,-3.250) -- (7.200,-4.150);
\draw[reduced,shift={(-0.047,-0.047)}] (0.951,-0.309) -- (0.150,0.300) -- (-1.176,1.618);
\draw[overlap edge,shift={(0.047,0.047)}] (0.951,-0.309) -- (0.150,0.300) -- (-1.176,1.618);
\draw[rounded corners=2pt,new edge] (0.588,0.809) -- (0.000,-1.000) -- (0.000,-1.900);
\draw[rounded corners=2pt,new edge] (9.012,0.809) -- (9.600,-1.000) -- (7.200,-3.250);
\draw[rounded corners=2pt,new edge] (8.649,-0.309) -- (9.450,0.300) -- (10.776,1.618);
\draw[rounded corners=2pt,new edge] (8.649,-0.309) -- (5.751,-0.309) -- (4.212,0.809);
\draw[rounded corners=2pt,new edge] (0.951,-0.309) -- (3.849,-0.309) -- (5.388,0.809);
\node[port,anchor=east] at (-2.015,-0.309) {$I_{0}$};
\node[port,anchor=east] at (-2.051,1.618) {$I_{1}$};
\node[port,anchor=north] at (0.000,-2.169) {$I_{2}$};
\node[port,anchor=west] at (11.663,1.618) {$O_{0}$};
\node[port,anchor=west] at (11.597,-0.309) {$O_{1}$};
\node[port,anchor=north] at (7.200,-4.419) {$O_{2}$};
\node[vertex dot] at (0.588,0.809) {};
\node[vertex dot] at (1.176,1.618) {};
\node[vertex dot] at (-0.951,-0.309) {};
\node[vertex dot] at (0.951,-0.309) {};
\node[vertex dot] at (-1.176,1.618) {};
\node[vertex dot] at (8.424,1.618) {};
\node[vertex dot] at (8.649,-0.309) {};
\node[vertex dot] at (5.976,1.618) {};
\node[vertex dot] at (5.388,0.809) {};
\node[vertex dot] at (4.212,0.809) {};
\node[vertex dot] at (4.800,-1.000) {};
\node[vertex dot] at (9.012,0.809) {};
\node[vertex dot] at (10.551,-0.309) {};
\node[vertex dot] at (10.776,1.618) {};
\node[vertex dot] at (3.624,1.618) {};
\node[vertex dot] at (7.200,-3.250) {};
\node[vlabel] at (0.153,1.143) {1};
\node[vlabel] at (1.516,2.060) {2};
\node[vlabel] at (-1.379,0.119) {3};
\node[vlabel] at (1.359,-0.717) {4};
\node[vlabel] at (-0.897,2.101) {5};
\node[vlabel] at (8.084,2.060) {8};
\node[vlabel] at (8.241,-0.717) {9};
\node[vlabel] at (6.356,2.114) {11};
\node[vlabel] at (5.983,0.562) {12};
\node[vlabel] at (3.550,0.720) {13};
\node[vlabel] at (4.710,-1.686) {14};
\node[vlabel] at (8.420,0.820) {16};
\node[vlabel] at (11.007,0.147) {17};
\node[vlabel] at (10.478,2.134) {18};
\node[vlabel] at (3.244,2.114) {19};
\node[vlabel] at (7.806,-3.250) {21};
\end{tikzpicture}
\end{minipage}
\hfill
\begin{minipage}[t]{0.49\textwidth}\centering
\begin{tikzpicture}[x=0.52cm,y=0.52cm,line cap=round,line join=round]
\tikzset{graph edge/.style={draw=black!48,line width=0.55pt},reduced/.style={draw=red!78!black,line width=1.8pt},new edge/.style={draw=green!50!black,line width=0.9pt,dash pattern=on 2.2pt off 1.35pt},overlap edge/.style={draw=green!50!black,line width=0.65pt,dash pattern=on 2.2pt off 1.35pt},vertex dot/.style={circle,fill=black,draw=black,inner sep=1.0pt,minimum size=2.0pt},vlabel/.style={font=\scriptsize,text=black!85,inner sep=0.3pt},port/.style={font=\scriptsize,inner sep=0.6pt}}
\node[font=\small,anchor=west] at (-2.450,3.050) {(c)};
\draw[graph edge] (0.588,0.809) -- (1.176,1.618);
\draw[graph edge] (0.588,0.809) -- (-0.951,-0.309);
\draw[reduced] (5.976,1.618) -- (5.388,0.809);
\draw[graph edge] (5.388,0.809) -- (4.800,-1.000);
\draw[graph edge] (4.212,0.809) -- (4.800,-1.000);
\draw[reduced] (8.649,-0.309) -- (10.551,-0.309);
\draw[graph edge] (9.012,0.809) -- (10.551,-0.309);
\draw[graph edge] (1.176,1.618) -- (3.624,1.618);
\draw[graph edge] (5.976,1.618) -- (3.624,1.618);
\draw[graph edge] (4.212,0.809) -- (3.624,1.618);
\draw[graph edge] (4.800,-1.000) -- (7.200,-3.250);
\draw[graph edge] (-0.951,-0.309) -- (-1.727,-0.309);
\draw[graph edge] (10.551,-0.309) -- (11.309,-0.309);
\draw[graph edge] (7.200,-3.250) -- (7.200,-4.150);
\draw[graph edge] (0.588,0.809) -- (0.000,-1.900);
\draw[graph edge] (9.012,0.809) -- (9.600,-1.000) -- (7.200,-3.250);
\draw[graph edge] (8.649,-0.309) -- (5.751,-0.309) -- (4.212,0.809);
\draw[rounded corners=2pt,new edge] (-0.951,-0.309) -- (4.200,-0.400) -- (5.388,0.809);
\draw[new edge] (1.176,1.618) -- (-1.176,1.618) -- (-1.763,1.618);
\draw[rounded corners=2pt,new edge] (5.976,1.618) -- (8.424,1.618) -- (9.012,0.809);
\draw[rounded corners=2pt,new edge] (8.649,-0.309) -- (10.776,1.618) -- (11.375,1.618);
\node[port,anchor=east] at (-2.015,-0.309) {$I_{0}$};
\node[port,anchor=east] at (-2.051,1.618) {$I_{1}$};
\node[port,anchor=north] at (0.000,-2.169) {$I_{2}$};
\node[port,anchor=west] at (11.663,1.618) {$O_{0}$};
\node[port,anchor=west] at (11.597,-0.309) {$O_{1}$};
\node[port,anchor=north] at (7.200,-4.419) {$O_{2}$};
\node[vertex dot] at (0.588,0.809) {};
\node[vertex dot] at (1.176,1.618) {};
\node[vertex dot] at (-0.951,-0.309) {};
\node[vertex dot] at (8.649,-0.309) {};
\node[vertex dot] at (5.976,1.618) {};
\node[vertex dot] at (5.388,0.809) {};
\node[vertex dot] at (4.212,0.809) {};
\node[vertex dot] at (4.800,-1.000) {};
\node[vertex dot] at (9.012,0.809) {};
\node[vertex dot] at (10.551,-0.309) {};
\node[vertex dot] at (3.624,1.618) {};
\node[vertex dot] at (7.200,-3.250) {};
\node[vlabel] at (0.153,1.143) {1};
\node[vlabel] at (1.516,2.060) {2};
\node[vlabel] at (-1.379,0.119) {3};
\node[vlabel] at (8.221,-0.737) {9};
\node[vlabel] at (6.356,2.114) {11};
\node[vlabel] at (5.954,0.482) {12};
\node[vlabel] at (3.550,0.720) {13};
\node[vlabel] at (4.710,-1.686) {14};
\node[vlabel] at (8.420,0.550) {16};
\node[vlabel] at (11.007,0.147) {17};
\node[vlabel] at (3.244,2.114) {19};
\node[vlabel] at (7.841,-3.422) {21};
\end{tikzpicture}
\end{minipage}
\par\medskip
\begin{minipage}[t]{0.32\textwidth}\centering
\begin{tikzpicture}[x=0.34cm,y=0.52cm,line cap=round,line join=round]
\tikzset{graph edge/.style={draw=black!48,line width=0.55pt},reduced/.style={draw=red!78!black,line width=1.8pt},new edge/.style={draw=green!50!black,line width=0.9pt,dash pattern=on 2.2pt off 1.35pt},overlap edge/.style={draw=green!50!black,line width=0.65pt,dash pattern=on 2.2pt off 1.35pt},vertex dot/.style={circle,fill=black,draw=black,inner sep=1.0pt,minimum size=2.0pt},vlabel/.style={font=\scriptsize,text=black!85,inner sep=0.3pt},port/.style={font=\scriptsize,inner sep=0.6pt}}
\node[font=\small,anchor=west] at (-2.450,3.050) {(d)};
\draw[graph edge] (0.588,0.809) -- (1.176,1.618);
\draw[graph edge] (0.588,0.809) -- (-0.951,-0.309);
\draw[graph edge] (4.212,0.809) -- (4.800,-1.000);
\draw[graph edge] (1.176,1.618) -- (3.624,1.618);
\draw[reduced] (4.212,0.809) -- (3.624,1.618);
\draw[graph edge] (4.800,-1.000) -- (7.200,-3.250);
\draw[graph edge] (-0.951,-0.309) -- (-1.727,-0.309);
\draw[graph edge] (7.200,-3.250) -- (7.200,-4.150);
\draw[graph edge] (0.588,0.809) -- (0.000,-1.900);
\draw[graph edge] (9.012,0.809) -- (9.600,-1.000) -- (7.200,-3.250);
\draw[graph edge] (1.176,1.618) -- (-1.763,1.618);
\draw[rounded corners=2pt,new edge] (4.212,0.809) -- (5.751,-0.309) -- (8.649,-0.309) -- (10.776,1.618) -- (11.375,1.618);
\draw[rounded corners=2pt,new edge] (9.012,0.809) -- (10.551,-0.309) -- (11.309,-0.309);
\draw[rounded corners=2pt,new edge] (9.012,0.809) -- (8.424,1.618) -- (3.624,1.618);
\draw[rounded corners=2pt,new edge] (-0.951,-0.309) -- (4.200,-0.400) -- (5.388,0.809) -- (4.800,-1.000);
\node[port,anchor=east] at (-2.168,-0.309) {$I_{0}$};
\node[port,anchor=east] at (-2.204,1.618) {$I_{1}$};
\node[port,anchor=north] at (0.000,-2.169) {$I_{2}$};
\node[port,anchor=west] at (11.816,1.618) {$O_{0}$};
\node[port,anchor=west] at (11.750,-0.309) {$O_{1}$};
\node[port,anchor=north] at (7.200,-4.419) {$O_{2}$};
\node[vertex dot] at (0.588,0.809) {};
\node[vertex dot] at (1.176,1.618) {};
\node[vertex dot] at (-0.951,-0.309) {};
\node[vertex dot] at (4.212,0.809) {};
\node[vertex dot] at (4.800,-1.000) {};
\node[vertex dot] at (9.012,0.809) {};
\node[vertex dot] at (3.624,1.618) {};
\node[vertex dot] at (7.200,-3.250) {};
\node[vlabel] at (-0.077,1.143) {1};
\node[vlabel] at (1.602,2.101) {2};
\node[vlabel] at (-1.470,0.133) {3};
\node[vlabel] at (3.450,0.650) {13};
\node[vlabel] at (4.406,-1.622) {14};
\node[vlabel] at (8.250,0.550) {16};
\node[vlabel] at (3.168,2.134) {19};
\node[vlabel] at (8.259,-3.250) {21};
\end{tikzpicture}
\end{minipage}
\hfill
\begin{minipage}[t]{0.32\textwidth}\centering
\begin{tikzpicture}[x=0.34cm,y=0.52cm,line cap=round,line join=round]
\tikzset{graph edge/.style={draw=black!48,line width=0.55pt},reduced/.style={draw=red!78!black,line width=1.8pt},new edge/.style={draw=green!50!black,line width=0.9pt,dash pattern=on 2.2pt off 1.35pt},overlap edge/.style={draw=green!50!black,line width=0.65pt,dash pattern=on 2.2pt off 1.35pt},vertex dot/.style={circle,fill=black,draw=black,inner sep=1.0pt,minimum size=2.0pt},vlabel/.style={font=\scriptsize,text=black!85,inner sep=0.3pt},port/.style={font=\scriptsize,inner sep=0.6pt}}
\node[font=\small,anchor=west] at (-2.450,3.050) {(e)};
\draw[graph edge] (0.588,0.809) -- (1.176,1.618);
\draw[graph edge] (0.588,0.809) -- (-0.951,-0.309);
\draw[graph edge] (4.800,-1.000) -- (7.200,-3.250);
\draw[graph edge] (-0.951,-0.309) -- (-1.727,-0.309);
\draw[graph edge] (7.200,-3.250) -- (7.200,-4.150);
\draw[graph edge] (0.588,0.809) -- (0.000,-1.900);
\draw[graph edge] (9.012,0.809) -- (9.600,-1.000) -- (7.200,-3.250);
\draw[graph edge] (1.176,1.618) -- (-1.763,1.618);
\draw[graph edge] (9.012,0.809) -- (11.309,-0.309);
\draw[graph edge] (-0.951,-0.309) -- (4.800,-1.000);
\draw[rounded corners=2pt,new edge] (4.800,-1.000) -- (4.212,0.809) -- (5.751,-0.309) -- (8.649,-0.309) -- (10.776,1.618) -- (11.375,1.618);
\draw[rounded corners=2pt,new edge] (1.176,1.618) -- (3.624,1.618) -- (8.424,1.618) -- (9.012,0.809);
\node[port,anchor=east] at (-2.168,-0.309) {$I_{0}$};
\node[port,anchor=east] at (-2.204,1.618) {$I_{1}$};
\node[port,anchor=north] at (0.000,-2.169) {$I_{2}$};
\node[port,anchor=west] at (11.816,1.618) {$O_{0}$};
\node[port,anchor=west] at (11.750,-0.309) {$O_{1}$};
\node[port,anchor=north] at (7.200,-4.419) {$O_{2}$};
\node[vertex dot] at (0.588,0.809) {};
\node[vertex dot] at (1.176,1.618) {};
\node[vertex dot] at (-0.951,-0.309) {};
\node[vertex dot] at (4.800,-1.000) {};
\node[vertex dot] at (9.012,0.809) {};
\node[vertex dot] at (7.200,-3.250) {};
\node[vlabel] at (-0.077,1.143) {1};
\node[vlabel] at (1.602,2.101) {2};
\node[vlabel] at (-1.470,0.133) {3};
\node[vlabel] at (4.406,-1.622) {14};
\node[vlabel] at (8.250,0.500) {16};
\node[vlabel] at (8.259,-3.250) {21};
\end{tikzpicture}
\end{minipage}
\hfill
\begin{minipage}[t]{0.32\textwidth}\centering
\begin{tikzpicture}[x=1cm,y=1cm,line cap=round,line join=round]
\tikzset{graph edge/.style={draw=black!48,line width=0.55pt},reduced/.style={draw=red!78!black,line width=1.8pt},new edge/.style={draw=green!50!black,line width=0.9pt,dash pattern=on 2.2pt off 1.35pt},overlap edge/.style={draw=green!50!black,line width=0.65pt,dash pattern=on 2.2pt off 1.35pt},vertex dot/.style={circle,fill=black,draw=black,inner sep=1.0pt,minimum size=2.0pt},vlabel/.style={font=\scriptsize,text=black!85,inner sep=0.3pt},port/.style={font=\scriptsize,inner sep=0.6pt}}
\node[font=\small,anchor=west] at (-1.950,1.586) {(f)};
\draw[graph edge] (-0.585,1.105) -- (0.585,1.105);
\draw[graph edge] (0.585,1.105) -- (0.585,-1.105);
\draw[graph edge] (0.585,-1.105) -- (0.234,-0.325);
\draw[new edge] (0.234,-0.325) -- (-0.234,-0.325);
\draw[graph edge] (-0.234,-0.325) -- (-0.585,-1.105);
\draw[graph edge] (-0.585,-1.105) -- (-0.585,1.105);
\draw[new edge] (0.585,1.105) -- (1.485,1.105);
\draw[graph edge] (0.585,-1.105) -- (1.485,-1.105);
\draw[graph edge] (0.234,-0.325) -- (0.559,0.000) -- (1.485,0.000);
\draw[graph edge] (-0.234,-0.325) -- (-0.559,0.000) -- (-1.485,0.000);
\draw[graph edge] (-0.585,-1.105) -- (-1.485,-1.105);
\draw[graph edge] (-0.585,1.105) -- (-1.485,1.105);
\node[port,anchor=east] at (-1.635,1.105) {$I_{0}$};
\node[port,anchor=west] at (1.635,1.105) {$O_{0}$};
\node[port,anchor=east] at (-1.635,-1.105) {$I_{2}$};
\node[port,anchor=west] at (1.635,-1.105) {$O_{2}$};
\node[port,anchor=east] at (-1.635,0.000) {$I_{1}$};
\node[port,anchor=west] at (1.635,0.000) {$O_{1}$};
\node[vertex dot] at (-0.585,1.105) {};
\node[vertex dot] at (0.585,1.105) {};
\node[vertex dot] at (0.585,-1.105) {};
\node[vertex dot] at (0.234,-0.325) {};
\node[vertex dot] at (-0.234,-0.325) {};
\node[vertex dot] at (-0.585,-1.105) {};
\node[vlabel] at (-0.808,-1.328) {1};
\node[vlabel] at (-0.194,-0.023) {2};
\node[vlabel] at (-0.808,1.328) {3};
\node[vlabel] at (0.822,1.342) {14};
\node[vlabel] at (0.130,0.039) {16};
\node[vlabel] at (0.822,-1.342) {21};
\end{tikzpicture}
\end{minipage}
\end{minipage}}
\par\vspace{1mm}
\textnormal{(g)}\par\vspace{1pt}
\resizebox{0.9\textwidth}{!}{\ringwiring{1}}
\caption{(a)--(e) Reductions of $H$ to a hexagon. Thick red edges are reduced
next; green dashed edges were created since the preceding panel. Thin green
dashes adjacent to a thick red edge indicate both roles. (f) The hexagon of (e),
redrawn. (g) As in Figure~\ref{fig:block}(c), with each block a hexagon in a
dashed box marking the former $H^i$.}
\label{fig:iext-block}
\end{figure}
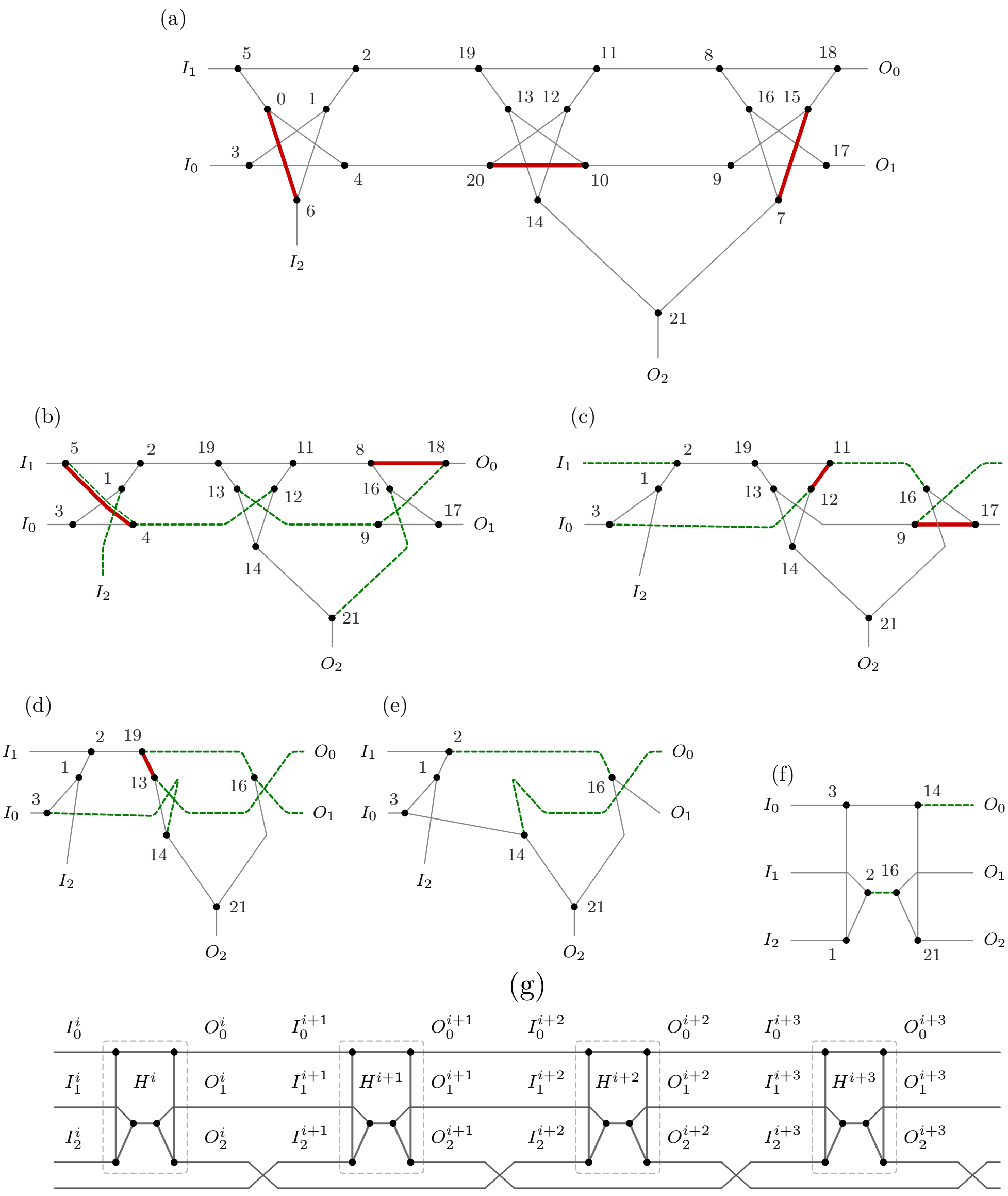

\begin{figure}[p]
\centering
\begin{minipage}{\textwidth}\centering\resizebox{0.97\textwidth}{!}{
\begin{tikzpicture}[x=1cm,y=1cm,line cap=round,line join=round]
  \tikzset{wire/.style={draw=black!65,line width=0.65pt},
    hexedge/.style={draw=black!58,line width=0.9pt},
    rededge/.style={draw=red!78!black,line width=1.8pt},
    greenedge/.style={draw=green!50!black,line width=0.9pt,dash pattern=on 2.2pt off 1.35pt},
    border/.style={draw=black!30,line width=0.5pt,dash pattern=on 2.4pt off 1.6pt,rounded corners=1.5pt},
    dotn/.style={circle,fill=black,inner sep=1.1pt},
    lab/.style={font=\footnotesize,inner sep=0pt},
    hlab/.style={font=\footnotesize,inner sep=0pt},
    plab/.style={font=\normalsize,inner sep=0pt}}
  \begin{scope}[yshift=0.00cm]
    \node[plab,anchor=east] at (-2.35,0.0) {(a)};
    \draw[border] (-0.650,-1.020) rectangle (0.650,0.670);
    \draw[border] (1.900,-1.020) rectangle (3.200,0.670);
    \draw[border] (4.450,-1.020) rectangle (5.750,0.670);
    \draw[border] (7.000,-1.020) rectangle (8.300,0.670);
    \draw[border] (9.550,-1.020) rectangle (10.850,0.670);
    \draw[rededge] (3.000,0.500) -- (3.000,-0.850);
    \draw[hexedge] (2.100,0.500) -- (3.000,0.500);
    \draw[wire] (3.000,0.500) -- (4.650,0.500);
    \draw[hexedge] (3.000,-0.850) -- (2.730,-0.250);
    \draw[wire] (3.000,-0.850) -- (3.600,-0.850) -- (4.050,-1.250) -- (6.150,-1.250) -- (6.600,-0.850) -- (7.200,-0.850);
    \draw[hexedge] (2.730,-0.250) -- (2.370,-0.250);
    \draw[wire] (2.730,-0.250) -- (2.980,0.000) -- (4.670,0.000) -- (4.920,-0.250);
    \draw[hexedge] (2.370,-0.250) -- (2.100,-0.850);
    \draw[wire] (0.180,-0.250) -- (0.430,0.000) -- (2.120,0.000) -- (2.370,-0.250);
    \draw[hexedge] (2.100,-0.850) -- (2.100,0.500);
    \draw[wire] (2.100,-0.850) -- (1.500,-0.850) -- (1.050,-1.250) -- (-1.050,-1.250) -- (-1.500,-0.850) -- (-1.750,-0.850);
    \draw[wire] (0.450,0.500) -- (2.100,0.500);
    \draw[hexedge] (5.550,0.500) -- (5.550,-0.850);
    \draw[hexedge] (4.650,0.500) -- (5.550,0.500);
    \draw[wire] (5.550,0.500) -- (7.200,0.500);
    \draw[hexedge] (5.550,-0.850) -- (5.280,-0.250);
    \draw[wire] (5.550,-0.850) -- (6.150,-0.850) -- (6.600,-1.250) -- (8.700,-1.250) -- (9.150,-0.850) -- (9.750,-0.850);
    \draw[hexedge] (5.280,-0.250) -- (4.920,-0.250);
    \draw[rededge] (5.280,-0.250) -- (5.530,0.000) -- (7.220,0.000) -- (7.470,-0.250);
    \draw[hexedge] (4.920,-0.250) -- (4.650,-0.850);
    \draw[hexedge] (4.650,-0.850) -- (4.650,0.500);
    \draw[wire] (0.450,-0.850) -- (1.050,-0.850) -- (1.500,-1.250) -- (3.600,-1.250) -- (4.050,-0.850) -- (4.650,-0.850);
    \draw[hexedge] (8.100,0.500) -- (8.100,-0.850);
    \draw[hexedge] (7.200,0.500) -- (8.100,0.500);
    \draw[wire] (8.100,0.500) -- (9.750,0.500);
    \draw[hexedge] (8.100,-0.850) -- (7.830,-0.250);
    \draw[wire] (11.950,-0.850) -- (11.700,-0.850) -- (11.250,-1.250) -- (9.150,-1.250) -- (8.700,-0.850) -- (8.100,-0.850);
    \draw[hexedge] (7.830,-0.250) -- (7.470,-0.250);
    \draw[wire] (7.830,-0.250) -- (8.080,0.000) -- (9.770,0.000) -- (10.020,-0.250);
    \draw[hexedge] (7.470,-0.250) -- (7.200,-0.850);
    \draw[hexedge] (7.200,-0.850) -- (7.200,0.500);
    \draw[hexedge] (10.650,0.500) -- (10.650,-0.850);
    \draw[hexedge] (9.750,0.500) -- (10.650,0.500);
    \draw[wire] (10.650,0.500) -- (11.950,0.500);
    \draw[hexedge] (10.650,-0.850) -- (10.380,-0.250);
    \draw[wire] (11.950,-1.250) -- (11.700,-1.250) -- (11.250,-0.850) -- (10.650,-0.850);
    \draw[hexedge] (10.380,-0.250) -- (10.020,-0.250);
    \draw[wire] (10.380,-0.250) -- (10.630,0.000) -- (11.200,0.000) -- (11.950,0.000);
    \draw[hexedge] (10.020,-0.250) -- (9.750,-0.850);
    \draw[hexedge] (9.750,-0.850) -- (9.750,0.500);
    \draw[wire] (-0.450,-0.850) -- (-1.050,-0.850) -- (-1.500,-1.250) -- (-1.750,-1.250);
    \draw[wire] (-1.750,0.500) -- (-0.450,0.500);
    \draw[wire] (-1.750,0.000) -- (-1.000,0.000) -- (-0.430,0.000) -- (-0.180,-0.250);
    \draw[hexedge] (0.450,0.500) -- (0.450,-0.850);
    \draw[hexedge] (-0.450,0.500) -- (0.450,0.500);
    \draw[hexedge] (0.450,-0.850) -- (0.180,-0.250);
    \draw[hexedge] (0.180,-0.250) -- (-0.180,-0.250);
    \draw[hexedge] (-0.180,-0.250) -- (-0.450,-0.850);
    \draw[hexedge] (-0.450,-0.850) -- (-0.450,0.500);
    \node[dotn] at (3.000,0.500) {};
    \node[dotn] at (3.000,-0.850) {};
    \node[dotn] at (2.730,-0.250) {};
    \node[dotn] at (2.370,-0.250) {};
    \node[dotn] at (2.100,-0.850) {};
    \node[dotn] at (2.100,0.500) {};
    \node[dotn] at (5.550,0.500) {};
    \node[dotn] at (5.550,-0.850) {};
    \node[dotn] at (5.280,-0.250) {};
    \node[dotn] at (4.920,-0.250) {};
    \node[dotn] at (4.650,-0.850) {};
    \node[dotn] at (4.650,0.500) {};
    \node[dotn] at (8.100,0.500) {};
    \node[dotn] at (8.100,-0.850) {};
    \node[dotn] at (7.830,-0.250) {};
    \node[dotn] at (7.470,-0.250) {};
    \node[dotn] at (7.200,-0.850) {};
    \node[dotn] at (7.200,0.500) {};
    \node[dotn] at (10.650,0.500) {};
    \node[dotn] at (10.650,-0.850) {};
    \node[dotn] at (10.380,-0.250) {};
    \node[dotn] at (10.020,-0.250) {};
    \node[dotn] at (9.750,-0.850) {};
    \node[dotn] at (9.750,0.500) {};
    \node[dotn] at (0.450,0.500) {};
    \node[dotn] at (0.450,-0.850) {};
    \node[dotn] at (0.180,-0.250) {};
    \node[dotn] at (-0.180,-0.250) {};
    \node[dotn] at (-0.450,-0.850) {};
    \node[dotn] at (-0.450,0.500) {};
    \node[hlab,anchor=south] at (0.000,0.840) {$H^{m-4}$};
    \node[hlab,anchor=south] at (2.550,0.840) {$H^{m-3}$};
    \node[hlab,anchor=south] at (5.100,0.840) {$H^{m-2}$};
    \node[hlab,anchor=south] at (7.650,0.840) {$H^{m-1}$};
    \node[hlab,anchor=south] at (10.200,0.840) {$H^{0}$};
  \end{scope}
  \begin{scope}[yshift=-4.33cm]
    \node[plab,anchor=east] at (-2.35,0.0) {(b)};
    \draw[border] (-0.650,-1.020) rectangle (0.650,0.670);
    \draw[border] (1.900,-1.020) rectangle (3.200,0.670);
    \draw[border] (4.450,-1.020) rectangle (5.750,0.670);
    \draw[border] (7.000,-1.020) rectangle (8.300,0.670);
    \draw[border] (9.550,-1.020) rectangle (10.850,0.670);
    \draw[hexedge] (2.730,-0.250) -- (2.370,-0.250);
    \draw[wire] (2.730,-0.250) -- (2.980,0.000) -- (4.670,0.000) -- (4.920,-0.250);
    \draw[greenedge] (2.730,-0.250) -- (3.000,-0.850) -- (3.600,-0.850) -- (4.050,-1.250) -- (6.150,-1.250) -- (6.600,-0.850) -- (7.200,-0.850);
    \draw[hexedge] (2.370,-0.250) -- (2.100,-0.850);
    \draw[wire] (0.180,-0.250) -- (0.430,0.000) -- (2.120,0.000) -- (2.370,-0.250);
    \draw[hexedge] (2.100,-0.850) -- (2.100,0.500);
    \draw[wire] (2.100,-0.850) -- (1.500,-0.850) -- (1.050,-1.250) -- (-1.050,-1.250) -- (-1.500,-0.850) -- (-1.750,-0.850);
    \draw[greenedge] (2.100,0.500) -- (3.000,0.500) -- (4.650,0.500);
    \draw[wire] (0.450,0.500) -- (2.100,0.500);
    \draw[hexedge] (5.550,0.500) -- (5.550,-0.850);
    \draw[hexedge] (4.650,0.500) -- (5.550,0.500);
    \draw[wire] (5.550,0.500) -- (7.200,0.500);
    \draw[greenedge] (4.920,-0.250) -- (5.280,-0.250) -- (5.550,-0.850);
    \draw[wire] (5.550,-0.850) -- (6.150,-0.850) -- (6.600,-1.250) -- (8.700,-1.250) -- (9.150,-0.850) -- (9.750,-0.850);
    \draw[hexedge] (4.920,-0.250) -- (4.650,-0.850);
    \draw[rededge] (4.650,-0.850) -- (4.650,0.500);
    \draw[wire] (0.450,-0.850) -- (1.050,-0.850) -- (1.500,-1.250) -- (3.600,-1.250) -- (4.050,-0.850) -- (4.650,-0.850);
    \draw[hexedge] (8.100,0.500) -- (8.100,-0.850);
    \draw[hexedge] (7.200,0.500) -- (8.100,0.500);
    \draw[wire] (8.100,0.500) -- (9.750,0.500);
    \draw[hexedge] (8.100,-0.850) -- (7.830,-0.250);
    \draw[wire] (11.950,-0.850) -- (11.700,-0.850) -- (11.250,-1.250) -- (9.150,-1.250) -- (8.700,-0.850) -- (8.100,-0.850);
    \draw[greenedge] (7.200,-0.850) -- (7.470,-0.250) -- (7.830,-0.250);
    \draw[wire] (7.830,-0.250) -- (8.080,0.000) -- (9.770,0.000) -- (10.020,-0.250);
    \draw[rededge] (7.200,-0.850) -- (7.200,0.500);
    \draw[hexedge] (10.650,0.500) -- (10.650,-0.850);
    \draw[hexedge] (9.750,0.500) -- (10.650,0.500);
    \draw[wire] (10.650,0.500) -- (11.950,0.500);
    \draw[hexedge] (10.650,-0.850) -- (10.380,-0.250);
    \draw[wire] (11.950,-1.250) -- (11.700,-1.250) -- (11.250,-0.850) -- (10.650,-0.850);
    \draw[hexedge] (10.380,-0.250) -- (10.020,-0.250);
    \draw[wire] (10.380,-0.250) -- (10.630,0.000) -- (11.200,0.000) -- (11.950,0.000);
    \draw[hexedge] (10.020,-0.250) -- (9.750,-0.850);
    \draw[hexedge] (9.750,-0.850) -- (9.750,0.500);
    \draw[wire] (-0.450,-0.850) -- (-1.050,-0.850) -- (-1.500,-1.250) -- (-1.750,-1.250);
    \draw[wire] (-1.750,0.500) -- (-0.450,0.500);
    \draw[wire] (-1.750,0.000) -- (-1.000,0.000) -- (-0.430,0.000) -- (-0.180,-0.250);
    \draw[hexedge] (0.450,0.500) -- (0.450,-0.850);
    \draw[hexedge] (-0.450,0.500) -- (0.450,0.500);
    \draw[hexedge] (0.450,-0.850) -- (0.180,-0.250);
    \draw[hexedge] (0.180,-0.250) -- (-0.180,-0.250);
    \draw[hexedge] (-0.180,-0.250) -- (-0.450,-0.850);
    \draw[hexedge] (-0.450,-0.850) -- (-0.450,0.500);
    \node[dotn] at (2.730,-0.250) {};
    \node[dotn] at (2.370,-0.250) {};
    \node[dotn] at (2.100,-0.850) {};
    \node[dotn] at (2.100,0.500) {};
    \node[dotn] at (5.550,0.500) {};
    \node[dotn] at (5.550,-0.850) {};
    \node[dotn] at (4.920,-0.250) {};
    \node[dotn] at (4.650,-0.850) {};
    \node[dotn] at (4.650,0.500) {};
    \node[dotn] at (8.100,0.500) {};
    \node[dotn] at (8.100,-0.850) {};
    \node[dotn] at (7.830,-0.250) {};
    \node[dotn] at (7.200,-0.850) {};
    \node[dotn] at (7.200,0.500) {};
    \node[dotn] at (10.650,0.500) {};
    \node[dotn] at (10.650,-0.850) {};
    \node[dotn] at (10.380,-0.250) {};
    \node[dotn] at (10.020,-0.250) {};
    \node[dotn] at (9.750,-0.850) {};
    \node[dotn] at (9.750,0.500) {};
    \node[dotn] at (0.450,0.500) {};
    \node[dotn] at (0.450,-0.850) {};
    \node[dotn] at (0.180,-0.250) {};
    \node[dotn] at (-0.180,-0.250) {};
    \node[dotn] at (-0.450,-0.850) {};
    \node[dotn] at (-0.450,0.500) {};
    \node[hlab,anchor=south] at (0.000,0.840) {$H^{m-4}$};
    \node[hlab,anchor=south] at (2.550,0.840) {$H^{m-3}$};
    \node[hlab,anchor=south] at (5.100,0.840) {$H^{m-2}$};
    \node[hlab,anchor=south] at (7.650,0.840) {$H^{m-1}$};
    \node[hlab,anchor=south] at (10.200,0.840) {$H^{0}$};
  \end{scope}
  \begin{scope}[yshift=-8.66cm]
    \node[plab,anchor=east] at (-2.35,0.0) {(c)};
    \draw[border] (-0.650,-1.020) rectangle (0.650,0.670);
    \draw[border] (1.900,-1.020) rectangle (3.200,0.670);
    \draw[border] (4.450,-1.020) rectangle (5.750,0.670);
    \draw[border] (7.000,-1.020) rectangle (8.300,0.670);
    \draw[border] (9.550,-1.020) rectangle (10.850,0.670);
    \draw[hexedge] (2.730,-0.250) -- (2.370,-0.250);
    \draw[rededge] (2.730,-0.250) -- (2.980,0.000) -- (4.670,0.000) -- (4.920,-0.250);
    \draw[greenedge] (2.730,-0.250) -- (3.000,-0.850) -- (3.600,-0.850) -- (4.050,-1.250) -- (6.150,-1.250) -- (6.600,-0.850) -- (7.200,-0.850) -- (7.470,-0.250) -- (7.830,-0.250);
    \draw[hexedge] (2.370,-0.250) -- (2.100,-0.850);
    \draw[wire] (0.180,-0.250) -- (0.430,0.000) -- (2.120,0.000) -- (2.370,-0.250);
    \draw[hexedge] (2.100,-0.850) -- (2.100,0.500);
    \draw[wire] (2.100,-0.850) -- (1.500,-0.850) -- (1.050,-1.250) -- (-1.050,-1.250) -- (-1.500,-0.850) -- (-1.750,-0.850);
    \draw[greenedge] (2.100,0.500) -- (3.000,0.500) -- (4.650,0.500) -- (5.550,0.500);
    \draw[wire] (0.450,0.500) -- (2.100,0.500);
    \draw[rededge] (5.550,0.500) -- (5.550,-0.850);
    \draw[greenedge] (5.550,0.500) -- (7.200,0.500) -- (8.100,0.500);
    \draw[hexedge] (4.920,-0.250) -- (5.280,-0.250) -- (5.550,-0.850);
    \draw[wire] (5.550,-0.850) -- (6.150,-0.850) -- (6.600,-1.250) -- (8.700,-1.250) -- (9.150,-0.850) -- (9.750,-0.850);
    \draw[greenedge] (0.450,-0.850) -- (1.050,-0.850) -- (1.500,-1.250) -- (3.600,-1.250) -- (4.050,-0.850) -- (4.650,-0.850) -- (4.920,-0.250);
    \draw[hexedge] (8.100,0.500) -- (8.100,-0.850);
    \draw[wire] (8.100,0.500) -- (9.750,0.500);
    \draw[hexedge] (8.100,-0.850) -- (7.830,-0.250);
    \draw[wire] (11.950,-0.850) -- (11.700,-0.850) -- (11.250,-1.250) -- (9.150,-1.250) -- (8.700,-0.850) -- (8.100,-0.850);
    \draw[wire] (7.830,-0.250) -- (8.080,0.000) -- (9.770,0.000) -- (10.020,-0.250);
    \draw[hexedge] (10.650,0.500) -- (10.650,-0.850);
    \draw[hexedge] (9.750,0.500) -- (10.650,0.500);
    \draw[wire] (10.650,0.500) -- (11.950,0.500);
    \draw[hexedge] (10.650,-0.850) -- (10.380,-0.250);
    \draw[wire] (11.950,-1.250) -- (11.700,-1.250) -- (11.250,-0.850) -- (10.650,-0.850);
    \draw[hexedge] (10.380,-0.250) -- (10.020,-0.250);
    \draw[wire] (10.380,-0.250) -- (10.630,0.000) -- (11.200,0.000) -- (11.950,0.000);
    \draw[hexedge] (10.020,-0.250) -- (9.750,-0.850);
    \draw[hexedge] (9.750,-0.850) -- (9.750,0.500);
    \draw[wire] (-0.450,-0.850) -- (-1.050,-0.850) -- (-1.500,-1.250) -- (-1.750,-1.250);
    \draw[wire] (-1.750,0.500) -- (-0.450,0.500);
    \draw[wire] (-1.750,0.000) -- (-1.000,0.000) -- (-0.430,0.000) -- (-0.180,-0.250);
    \draw[hexedge] (0.450,0.500) -- (0.450,-0.850);
    \draw[hexedge] (-0.450,0.500) -- (0.450,0.500);
    \draw[hexedge] (0.450,-0.850) -- (0.180,-0.250);
    \draw[hexedge] (0.180,-0.250) -- (-0.180,-0.250);
    \draw[hexedge] (-0.180,-0.250) -- (-0.450,-0.850);
    \draw[hexedge] (-0.450,-0.850) -- (-0.450,0.500);
    \node[dotn] at (2.730,-0.250) {};
    \node[dotn] at (2.370,-0.250) {};
    \node[dotn] at (2.100,-0.850) {};
    \node[dotn] at (2.100,0.500) {};
    \node[dotn] at (5.550,0.500) {};
    \node[dotn] at (5.550,-0.850) {};
    \node[dotn] at (4.920,-0.250) {};
    \node[dotn] at (8.100,0.500) {};
    \node[dotn] at (8.100,-0.850) {};
    \node[dotn] at (7.830,-0.250) {};
    \node[dotn] at (10.650,0.500) {};
    \node[dotn] at (10.650,-0.850) {};
    \node[dotn] at (10.380,-0.250) {};
    \node[dotn] at (10.020,-0.250) {};
    \node[dotn] at (9.750,-0.850) {};
    \node[dotn] at (9.750,0.500) {};
    \node[dotn] at (0.450,0.500) {};
    \node[dotn] at (0.450,-0.850) {};
    \node[dotn] at (0.180,-0.250) {};
    \node[dotn] at (-0.180,-0.250) {};
    \node[dotn] at (-0.450,-0.850) {};
    \node[dotn] at (-0.450,0.500) {};
    \node[hlab,anchor=south] at (0.000,0.840) {$H^{m-4}$};
    \node[hlab,anchor=south] at (2.550,0.840) {$H^{m-3}$};
    \node[hlab,anchor=south] at (5.100,0.840) {$H^{m-2}$};
    \node[hlab,anchor=south] at (7.650,0.840) {$H^{m-1}$};
    \node[hlab,anchor=south] at (10.200,0.840) {$H^{0}$};
  \end{scope}
  \begin{scope}[yshift=-12.99cm]
    \node[plab,anchor=east] at (-2.35,0.0) {(d)};
    \draw[border] (-0.650,-1.020) rectangle (0.650,0.670);
    \draw[border] (9.550,-1.020) rectangle (10.850,0.670);
    \draw[hexedge] (2.370,-0.250) -- (2.100,-0.850);
    \draw[greenedge] (7.830,-0.250) -- (7.470,-0.250) -- (7.200,-0.850) -- (6.600,-0.850) -- (6.150,-1.250) -- (4.050,-1.250) -- (3.600,-0.850) -- (3.000,-0.850) -- (2.730,-0.250) -- (2.370,-0.250);
    \draw[wire] (0.180,-0.250) -- (0.430,0.000) -- (2.120,0.000) -- (2.370,-0.250);
    \draw[hexedge] (2.100,-0.850) -- (2.100,0.500);
    \draw[wire] (2.100,-0.850) -- (1.500,-0.850) -- (1.050,-1.250) -- (-1.050,-1.250) -- (-1.500,-0.850) -- (-1.750,-0.850);
    \draw[greenedge] (2.100,0.500) -- (3.000,0.500) -- (4.650,0.500) -- (5.550,0.500) -- (7.200,0.500) -- (8.100,0.500);
    \draw[wire] (0.450,0.500) -- (2.100,0.500);
    \draw[hexedge] (8.100,0.500) -- (8.100,-0.850);
    \draw[wire] (8.100,0.500) -- (9.750,0.500);
    \draw[hexedge] (8.100,-0.850) -- (7.830,-0.250);
    \draw[wire] (11.950,-0.850) -- (11.700,-0.850) -- (11.250,-1.250) -- (9.150,-1.250) -- (8.700,-0.850) -- (8.100,-0.850);
    \draw[wire] (7.830,-0.250) -- (8.080,0.000) -- (9.770,0.000) -- (10.020,-0.250);
    \draw[hexedge] (10.650,0.500) -- (10.650,-0.850);
    \draw[hexedge] (9.750,0.500) -- (10.650,0.500);
    \draw[wire] (10.650,0.500) -- (11.950,0.500);
    \draw[hexedge] (10.650,-0.850) -- (10.380,-0.250);
    \draw[wire] (11.950,-1.250) -- (11.700,-1.250) -- (11.250,-0.850) -- (10.650,-0.850);
    \draw[hexedge] (10.380,-0.250) -- (10.020,-0.250);
    \draw[wire] (10.380,-0.250) -- (10.630,0.000) -- (11.200,0.000) -- (11.950,0.000);
    \draw[hexedge] (10.020,-0.250) -- (9.750,-0.850);
    \draw[hexedge] (9.750,-0.850) -- (9.750,0.500);
    \draw[greenedge] (9.750,-0.850) -- (9.150,-0.850) -- (8.700,-1.250) -- (6.600,-1.250) -- (6.150,-0.850) -- (5.550,-0.850) -- (5.280,-0.250) -- (4.920,-0.250) -- (4.650,-0.850) -- (4.050,-0.850) -- (3.600,-1.250) -- (1.500,-1.250) -- (1.050,-0.850) -- (0.450,-0.850);
    \draw[wire] (-0.450,-0.850) -- (-1.050,-0.850) -- (-1.500,-1.250) -- (-1.750,-1.250);
    \draw[wire] (-1.750,0.500) -- (-0.450,0.500);
    \draw[wire] (-1.750,0.000) -- (-1.000,0.000) -- (-0.430,0.000) -- (-0.180,-0.250);
    \draw[hexedge] (0.450,0.500) -- (0.450,-0.850);
    \draw[hexedge] (-0.450,0.500) -- (0.450,0.500);
    \draw[hexedge] (0.450,-0.850) -- (0.180,-0.250);
    \draw[hexedge] (0.180,-0.250) -- (-0.180,-0.250);
    \draw[hexedge] (-0.180,-0.250) -- (-0.450,-0.850);
    \draw[hexedge] (-0.450,-0.850) -- (-0.450,0.500);
    \node[dotn] at (2.370,-0.250) {};
    \node[dotn] at (2.100,-0.850) {};
    \node[dotn] at (2.100,0.500) {};
    \node[dotn] at (8.100,0.500) {};
    \node[dotn] at (8.100,-0.850) {};
    \node[dotn] at (7.830,-0.250) {};
    \node[dotn] at (10.650,0.500) {};
    \node[dotn] at (10.650,-0.850) {};
    \node[dotn] at (10.380,-0.250) {};
    \node[dotn] at (10.020,-0.250) {};
    \node[dotn] at (9.750,-0.850) {};
    \node[dotn] at (9.750,0.500) {};
    \node[dotn] at (0.450,0.500) {};
    \node[dotn] at (0.450,-0.850) {};
    \node[dotn] at (0.180,-0.250) {};
    \node[dotn] at (-0.180,-0.250) {};
    \node[dotn] at (-0.450,-0.850) {};
    \node[dotn] at (-0.450,0.500) {};
    \node[hlab,anchor=south] at (0.000,0.840) {$H^{m-4}$};
    \node[hlab,anchor=south] at (10.200,0.840) {$H^{0}$};
  \end{scope}
  \begin{scope}[yshift=-17.32cm]
    \node[plab,anchor=east] at (-2.35,0.0) {(e)};
    \draw[border] (1.900,-1.020) rectangle (8.300,0.670);
    \draw[border] (-0.650,-1.020) rectangle (0.650,0.670);
    \draw[border] (9.550,-1.020) rectangle (10.850,0.670);
    \draw[hexedge] (2.370,-0.250) -- (2.100,-0.850);
    \draw[greenedge] (7.830,-0.250) -- (2.370,-0.250);
    \draw[wire] (0.180,-0.250) -- (0.430,0.000) -- (2.120,0.000) -- (2.370,-0.250);
    \draw[hexedge] (2.100,-0.850) -- (2.100,0.500);
    \draw[wire] (2.100,-0.850) -- (1.500,-0.850) -- (1.050,-1.250) -- (-1.050,-1.250) -- (-1.500,-0.850) -- (-1.750,-0.850);
    \draw[greenedge] (2.100,0.500) -- (8.100,0.500);
    \draw[wire] (0.450,0.500) -- (2.100,0.500);
    \draw[hexedge] (8.100,0.500) -- (8.100,-0.850);
    \draw[wire] (8.100,0.500) -- (9.750,0.500);
    \draw[hexedge] (8.100,-0.850) -- (7.830,-0.250);
    \draw[wire] (11.950,-0.850) -- (11.700,-0.850) -- (11.250,-1.250) -- (9.150,-1.250) -- (8.700,-0.850) -- (8.100,-0.850);
    \draw[wire] (7.830,-0.250) -- (8.080,0.000) -- (9.770,0.000) -- (10.020,-0.250);
    \draw[hexedge] (10.650,0.500) -- (10.650,-0.850);
    \draw[hexedge] (9.750,0.500) -- (10.650,0.500);
    \draw[wire] (10.650,0.500) -- (11.950,0.500);
    \draw[hexedge] (10.650,-0.850) -- (10.380,-0.250);
    \draw[wire] (11.950,-1.250) -- (11.700,-1.250) -- (11.250,-0.850) -- (10.650,-0.850);
    \draw[hexedge] (10.380,-0.250) -- (10.020,-0.250);
    \draw[wire] (10.380,-0.250) -- (10.630,0.000) -- (11.200,0.000) -- (11.950,0.000);
    \draw[hexedge] (10.020,-0.250) -- (9.750,-0.850);
    \draw[hexedge] (9.750,-0.850) -- (9.750,0.500);
    \draw[greenedge] (9.750,-0.850) -- (9.150,-0.850) -- (8.700,-1.250) -- (1.500,-1.250) -- (1.050,-0.850) -- (0.450,-0.850);
    \draw[wire] (-0.450,-0.850) -- (-1.050,-0.850) -- (-1.500,-1.250) -- (-1.750,-1.250);
    \draw[wire] (-1.750,0.500) -- (-0.450,0.500);
    \draw[wire] (-1.750,0.000) -- (-1.000,0.000) -- (-0.430,0.000) -- (-0.180,-0.250);
    \draw[hexedge] (0.450,0.500) -- (0.450,-0.850);
    \draw[hexedge] (-0.450,0.500) -- (0.450,0.500);
    \draw[hexedge] (0.450,-0.850) -- (0.180,-0.250);
    \draw[hexedge] (0.180,-0.250) -- (-0.180,-0.250);
    \draw[hexedge] (-0.180,-0.250) -- (-0.450,-0.850);
    \draw[hexedge] (-0.450,-0.850) -- (-0.450,0.500);
    \node[dotn] at (2.370,-0.250) {};
    \node[dotn] at (2.100,-0.850) {};
    \node[dotn] at (2.100,0.500) {};
    \node[dotn] at (8.100,0.500) {};
    \node[dotn] at (8.100,-0.850) {};
    \node[dotn] at (7.830,-0.250) {};
    \node[dotn] at (10.650,0.500) {};
    \node[dotn] at (10.650,-0.850) {};
    \node[dotn] at (10.380,-0.250) {};
    \node[dotn] at (10.020,-0.250) {};
    \node[dotn] at (9.750,-0.850) {};
    \node[dotn] at (9.750,0.500) {};
    \node[dotn] at (0.450,0.500) {};
    \node[dotn] at (0.450,-0.850) {};
    \node[dotn] at (0.180,-0.250) {};
    \node[dotn] at (-0.180,-0.250) {};
    \node[dotn] at (-0.450,-0.850) {};
    \node[dotn] at (-0.450,0.500) {};
    \node[hlab,anchor=south] at (0.000,0.840) {$H^{m-4}$};
    \node[hlab,anchor=south] at (10.200,0.840) {$H^{0}$};
    \node[hlab,anchor=south] at (5.100,0.840) {new $H^{m-3}$};
  \end{scope}
\end{tikzpicture}
}\end{minipage}
\caption{The step from $R_m$ to $R_{m-2}$ on blocks $m-3,m-2,m-1$, drawn as in
Figure~\ref{fig:iext-block}(g). Thick red edges are reduced next; green dashed
edges were created since the preceding panel, drawn along the old route through
the deleted vertices. Blocks $m-4$ and $0$ are unchanged. In (e)
the green edges of (d) are straightened, and the three blocks have become the new
hexagon $H^{m-3}$.}
\label{fig:iext-ring-step}
\end{figure}
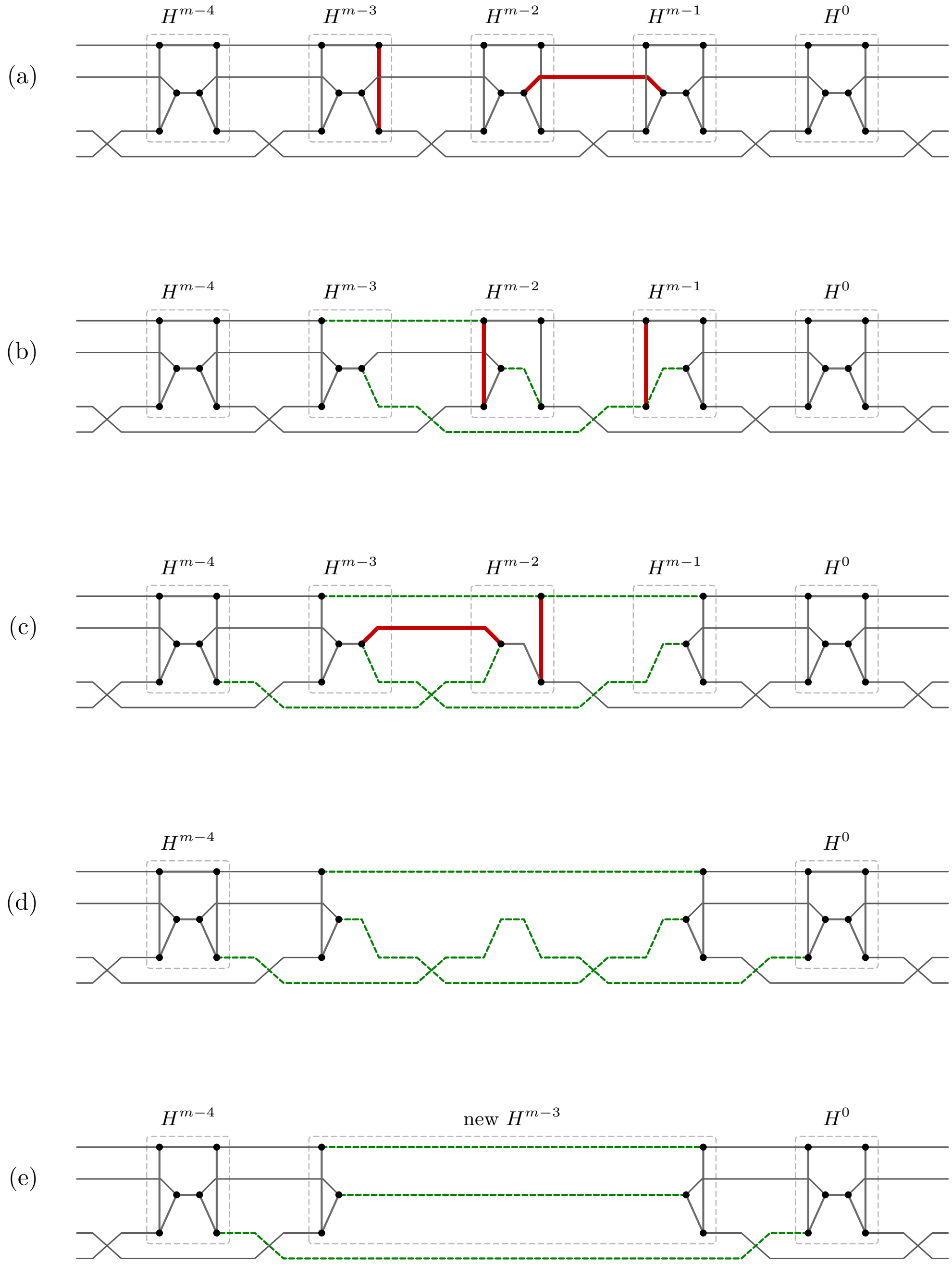

For $m\ge5$, the reductions in Figure~\ref{fig:iext-ring-step} are all valid
and take $R_m$ to $R_{m-2}$: blocks $m-3,m-2,m-1$ merge into a new block $H^{m-3}$, a
hexagon whose $I_j,O_j$ are the former $I_j^{m-3},O_j^{m-1}$ for
$j\in\{0,1,2\}$; every other block keeps its index.

The graphs $R_3$ and $R_4$ are cyclically $5$-connected (checked by a computer:
both reduce to the Petersen graph by a series of valid I-reductions, but we also
computed the cyclic connectivity directly). Applying Lemma~\ref{lem:iext} along the
reversed reduction sequences shows that $G_m$ is cyclically $5$-connected
for every $m\ge3$.
\end{proof}

\clearpage
\section*{Acknowledgements}

This work was partially supported from the research grants APVV-23-0076, VEGA 1/0173/25, and
VEGA 1/0613/26.
The author used large language models from Anthropic and OpenAI to assist with writing the manuscript and with developing proofs and code, and takes full responsibility for its content.


\begin{thebibliography}{99}
\bibitem{BJMY}
M. Bilinski, B. Jackson, J. Ma, and X. Yu, \emph{Circumference of $3$-connected claw-free graphs and large Eulerian subgraphs of $3$-edge-connected graphs}, J. Combin. Theory Ser. B \textbf{101} (2011), 214--236.\\
\doi{10.1016/j.jctb.2011.02.009}
\bibitem{BS}
J. A. Bondy and M. Simonovits, \emph{Longest cycles in $3$-connected cubic graphs}, Canad. J. Math. \textbf{32} (1980), 987--992.\\
\doi{10.4153/CJM-1980-076-2}
\bibitem{J}
B. Jackson, \emph{Longest cycles in $3$-connected cubic graphs}, J. Combin. Theory Ser. B \textbf{41} (1986), 17--26.\\
\doi{10.1016/0095-8956(86)90024-9}
\bibitem{LSVZ}
O.-H. S. Lo, J. M. Schmidt, N. Van Cleemput, and C. T. Zamfirescu, \emph{Shortness coefficient of cyclically $4$-edge-connected cubic graphs}, Electron. J. Combin. \textbf{27} (2020), \#P1.43.
\bibitem{LYZ}
Q. Liu, X. Yu, and Z. Zhang, \emph{Circumference of $3$-connected cubic graphs}, J. Combin. Theory Ser. B \textbf{128} (2018), 134--159.\\
\doi{10.1016/j.jctb.2017.08.008}
\bibitem{MM}
E. M\'a\v{c}ajov\'a and J. Maz\'ak, \emph{Cubic graphs with large circumference deficit}, J. Graph Theory \textbf{82} (2016), 433--440.\\
\doi{10.1002/jgt.21911}
\bibitem{MMpm}
E. M\'a\v{c}ajov\'a and J. Maz\'ak, \emph{On covering cubic graphs with three perfect matchings}, Discuss. Math. Graph Theory \textbf{46} (2026), 599--608.\\
\doi{10.7151/dmgt.2622}
\bibitem{MS}
E. M\'a\v{c}ajov\'a and M. \v{S}koviera, \emph{Irreducible snarks of given order and cyclic connectivity}, Discrete Math. \textbf{306} (2006), 779--791.\\
\doi{10.1016/j.disc.2006.02.003}
\bibitem{MRS}
J. Maz\'ak, J. Rajn\'ik, and M. \v{S}koviera,
\emph{Morphology of small snarks},
Electron. J. Combin. {\bf 29}(4) (2022), \#P4.30.\\
\doi{10.37236/10917}
\end{thebibliography}
\end{document}